\documentclass[12pt]{article}
\usepackage{unformatted}
\usepackage{url}
\date{September 1, 2020}
\usepackage{hyperref}
\hypersetup{
    colorlinks=true,
    linkcolor=blue,
    filecolor=magenta,
    urlcolor=cyan,
    citecolor=blue,
}
\usepackage[numbers]{natbib}
\usepackage{a4}
\usepackage{ifthen}
\renewcommand{\cite}[2][nobrackets]{\ifthenelse{\equal{#1}{nobrackets}}{\Citet{#2}}{\Citep{#2}}}

\newcommand{\articletype}[1]{}
\newcommand{\received}[1]{}
\newcommand{\revised}[1]{}
\newcommand{\accepted}[1]{}
\newcommand{\corres}[1]{}
\newcommand{\orgdiv}[1]{#1}
\newcommand{\orgname}[1]{#1}
\newcommand{\orgaddress}[1]{#1}
\newcommand{\country}[1]{#1}

\articletype{Research article}%

\received{}
\revised{}
\accepted{}

\usepackage{amsmath}
\usepackage{amssymb}
\usepackage{amsthm}

\usepackage{todonotes}

\usepackage{graphicx}

\usepackage{tikz}
\usepackage{pgfplots}
\usepackage{xcolor}

\usepackage{url}
\usepackage{datetime}

\usepackage{subfig}
\let\subtop\subfloat

\newcommand{\pd}[2]{\frac{\partial #1}{\partial #2}}
\newcommand{\bvec}[1]{\textbf{#1}}
\def\txtd{{\operatorname{d}}}

\renewcommand{\eqref}[1]{(\ref{eq:#1})}
\newcommand{\figref}[1]{Figure \ref{fig:#1}}
\newcommand{\tabref}[1]{Table \ref{tab:#1}}
\newcommand{\secref}[1]{Section \ref{sec:#1}}

\newlength{\figureWidth}
\newlength{\figureHeight}

\def\R{\mathbb{R}}

\begin{document}

\title{A staggered-grid multilevel incomplete LU for steady incompressible flows}

\address[rug]{\orgdiv{Bernoulli Institute for Mathematics, Computer Science and Artificial Intelligence}, \orgname{University of Groningen}, \orgaddress{Groningen, \country{the Netherlands}}}
\address[rug2]{\orgdiv{Zernike Institute for Advanced Materials}, \orgname{University of Groningen}, \orgaddress{Groningen, \country{the Netherlands}}}
\address[dlr]{\orgname{German Aerospace Center (DLR)}, \orgdiv{Simulation and Software Technology}, \orgaddress{Linder H\"ohe, 51147 Cologne, \country{Germany}}}

\author[rug]{Sven Baars}\ead{s.baars@rug.nl}
\author[rug2]{Mark van der Klok}
\author[dlr]{Jonas Thies}
\author[rug]{Fred W. Wubs}

\corres{Sven Baars, Bernoulli Institute for Mathematics, Computer Science and Artificial Intelligence,
University of Groningen,
Nijenborgh 9, 9747 AG Groningen,
The Netherlands,
\ead{s.baars@rug.nl}}

\begin{abstract}
Algorithms for studying transitions and instabilities in incompressible flows typically require the solution of linear systems with the full Jacobian matrix. 
Other popular approaches, like gradient-based design optimization and fully implicit time integration, also require very robust solvers for this type of linear system.
We present a parallel fully coupled multilevel incomplete factorization preconditioner for the 3D stationary incompressible Navier--Stokes equations on a structured grid.
The algorithm and software are based on the robust two-level method developed in \citet{wubs:11}.
In this paper, we identify some of the weak spots of the two-level scheme and propose remedies such as a different domain partitioning and recursive application of the method.
We apply the method to the well-known 3D lid-driven cavity benchmark problem, and demonstrate its superior robustness by comparing with a segregated SIMPLE-type preconditioner.\end{abstract}


\maketitle

\section{Introduction}
The incompressible Navier--Stokes equations accurately describe flow of Newtonian fluids like water and air at low Mach numbers.
In this paper we consider the numerical solution of these equations after spatial discretization on a structured grid.
The discretization we use does not introduce artificial diffusion and is therefore popular for direct numerical simulations of turbulent flow.
In this flow regime, one is typically interested in accurately resolving the temporal evolution of the flow.
The relatively small time steps in such a simulation allow simplifications in the solution process, in particular, Picard linearization (instead of Newton's), and a segregated solution scheme that treats the velocities and pressure separately \cite{verstappen:03}.

On the other end of the spectrum are Stokes flows and flows at very low Reynolds numbers.
Here it is common to solve the coupled linear systems with a segregated preconditioner, see e.g.\ \citet{elman:14}.
This class of linear solvers tries to reduce the problem to scalar linear systems that can be solved by multigrid methods.

In this paper, we are interested in the intermediate flow regimes, where the focus lies on accurately desribing flow transitions occuring when model parameters are varied (or uncertain).
In such situations it may be beneficial to directly compute steady state solutions, or at least take very large implicit time steps in order to quickly reach a statistical equilibrium solution.
Another application where the solution of linear systems with the full Jacobian (and its adjoint) play a key role is design optimization with gradient-based optimization methods \cite{brandenburg:09}.
In that application, however, unstructured grids may be preferable, which we do not address here.

To study the stability of steady states in fluid flow problems as a function of a parameter, e.g.\ the Reynolds number (see \citet{charru:11} for an introduction to such problems), a continuation method \citep{keller:77} can be used.
To apply such a method to high-dimensional systems, we require an efficient and robust method for solving linear systems associated with the discretized incompressible Navier--Stokes equations.
An elegant way of solving these systems is by replacing the complete LU factorization by an accurate incomplete one, which can be used as a preconditioner in an iterative procedure.
By an appropriate ordering technique and dropping procedure, one can construct an incomplete LU (ILU) factorization that yields grid independent convergence behavior and approaches an exact factorization as the amount of allowed fill is increased.
During the continuation process, this preconditioner can also be used to find approximate smallest eigenvalues and eigenvectors of the Jacobian matrix of the incompressible Navier--Stokes equations \citet{sleijpen:04}.
These eigenvalues are of interest since a switch in sign may indicate a bifurcation.

\begin{figure}[t]
  \begin{center}
    \begin{tikzpicture}
      \node (img1) {\includegraphics[width=.15\textwidth]{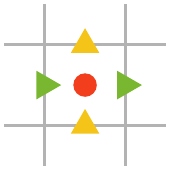}};
      \node (u) at (img1.west) [xshift=.5cm] {$u$};
      \node (v) at (img1.north) [yshift=-.4cm] {$v$};
      \node (p) at (img1) [xshift=-.2cm, yshift=.2cm] {$p$};
    \end{tikzpicture}
  \end{center}
  \caption{Positioning of unknowns in an Arakawa C-grid.} \label{fig:cgrid}
\end{figure}
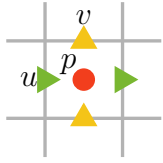

The incompressible Navier--Stokes equations can be written as
\begin{align}\label{eq:smilu-ns}
\begin{split}
  \pd{\bvec{u}}{t}+\bvec{u}\cdot\nabla\bvec{u} &= -\nabla p+\frac{1}{\text{Re}}\Delta \bvec{u},\\
  \nabla \cdot \bvec{u} &= 0,
  \end{split}
\end{align}
where $\text{Re}=\frac{\rho UD}{\mu}$ is the Reynolds number, $\rho$ is the density and $\mu$ is the dynamic viscosity, and $D$ and $U$ are characteristic length and velocity scales of the flow.
These equations are discretized using a second-order symmetry-preserving finite volume method on an Arakawa 
C-grid \cite{arakawa:77}; see \figref{cgrid}.
The discretization leads to a system of ordinary differential equations (ODEs)
\begin{align*}
  M\frac{\txtd \bvec{u}}{\txtd t} + N(\bvec{u},\bvec{u}) &= - G\bvec{p} + \frac{1}{\text{Re}} L \bvec{u}  + \bvec{f}_{\bvec{u}},\\
  -G^T \bvec{u} &= \bvec{f}_p,
\end{align*}
where $\bvec{u}$ and $\bvec{p}$ now represent the velocity and pressure in each grid point, $N(\cdot,\cdot)$ is the bilinear form arising from the convective terms, $L$ is the discretized Laplace operator, $G$ is the discretized gradient operator, $M$ is the mass matrix, which contains the volumes of the grid cells on its diagonal, and $\bvec{f}$ contains the known part of the boundary conditions.
Our method will exploit the property that the divergence operator is given by $-G^T$.

The term $N(\bvec{u},\bvec{v})$ is the discretized form of $\bvec{u} \cdot \nabla \bvec{v}$, and is bilinear. 
Hence for given $\bvec{u}$ the expression is linear in $\bvec{v}$: there exists a matrix $N_1(\bvec{u})$ such that $N(\bvec{u},\bvec{v})=N_1(\bvec{u})\bvec{v}$. 
Similarly, for given $\bvec{v}$, there exists a  matrix $N_2(\bvec{v})$ such that $N(\bvec{u},\bvec{v})=N_2(\bvec{v})\bvec{u}$, which is the discretized form of $\bvec{v} (\nabla \cdot \bvec{u})$. 
For the contribution of $N(\bvec{u},\bvec{v})$ to the Jacobian, we consider the derivative of  $N(\bvec{u},\bvec{u})$ in the direction $\Delta \bvec{u}$, which is found from the following expression by taking the limit $\epsilon\rightarrow 0$:
\begin{align*}
&[N(\bvec{u}+\epsilon\Delta \bvec{u},\bvec{u}+ \epsilon \Delta \bvec{u}) -N(\bvec{u},\bvec{u})]/\epsilon \\
&=N(\bvec{u},\Delta \bvec{u})+  N(\Delta \bvec{u},\bvec{u}) + \epsilon N(\Delta \bvec{u},\Delta \bvec{u}) \\
&= N_1(\bvec{u})\Delta \bvec{u}+  N_2(\bvec{u})\Delta \bvec{u} + \epsilon N(\Delta \bvec{u},\Delta \bvec{u}),
\end{align*}
where we used only the bilinearity of the expression $N(\cdot,\cdot)$. The last term becomes zero when taking the limit.

Using this notation, the linear system of saddle point type \citep{benzi:05} that has to be solved within each Newton step is given by
\begin{align}\label{eq:ns_spp}
  \begin{pmatrix}
    N_1(\bvec{u}) + N_2(\bvec{u}) - \frac{1}{\text{Re}} L & G \\
    G^T & 0
  \end{pmatrix}
  \begin{pmatrix}
    \Delta \bvec{u} \\
    \Delta \bvec{p}
  \end{pmatrix}
  = -
  \begin{pmatrix}
    \bvec{f}_{\bvec{u}} \\
    \bvec{f}_p
  \end{pmatrix}.
\end{align}

\setlength{\figureWidth}{0.3\textwidth}
\setlength{\figureHeight}{0.3\textwidth}

It is known \citep{verstappen:03} that, on closed domains, dissipation of kinetic energy only occurs by diffusion.
Our discretization preserves this property, which has the consequence that for divergence-free $\bvec{u}$, the operator $N_1(\bvec{u})$ is skew-symmetric.
A popular simplification is to omit $N_2(\bvec{u})$, which replaces the Newton process by a Picard iteration 
\citep{elman:14}. The approximate Jacobian then becomes negative definite (on the space of discrete divergence-free 
velocities), which greatly simplifies the solution process since definiteness is a condition under which many standard iterative methods converge.
We remark, however, that Picard iteration works well only for fairly low Reynolds numbers and far away from bifurcation points where steady solutions become unstable,
and seriously impairs the convergence rate of the nonlinear iteration \citep{carey:87, elman:03}.
Furthermore, some authors use time dependent approaches to study the stability of steady states \citep{dijkstra:14}.
This approach, however, also requires some tricks to obtain the desired information.
Since we want to study precisely the phenomena where the above methods experience difficulties, we would rather use the full Jacobian matrix of the nonlinear equations, applying Newton's method.

There are many methods, based on segregation of physical variables, that can solve the linear equations that arise in every Newton iteration.
In this approach the system is split into subproblems of an easier form for which standard methods exist.
The segregation can already be done at the discretization level, e.g.\ by doing a time integration and solving a pressure-correction equation independently of the momentum equations \citep{verstappen:03, verstappen:05}.
Another class of methods performs the segregation during the linear system solve, often in a preconditioning step.
Physics based preconditioners \citep{niet:07, klaij:12, cyr:12, he:18} try to split the problem into subsystems which capture the bulk of the physics.
The subsystems are again solved by iterative procedures, e.g.\ algebraic multigrid (AMG) for Poisson-like equations.
These methods consist of nested loops for: the nonlinear iteration, iterations over the coupled system, and iterations over the subsystems.
The stopping criteria of all these different iterations have to be tuned to make the solver efficient.
Furthermore, the total number of iterations in the innermost loop is given by the product of the number of iterations performed on all three levels of iteration and thus easily becomes excessive.
This is a major problem when trying to achieve extreme parallelism, because each innermost iteration typically requires global communication in the inner products.
The number of levels of nested iteration may increase even more if additional physical phenomena are added \citep{niet:07, thies:09}.
Our ultimate aim is to remove the inner iterations altogether and to solve the nonlinear equations using a subspace accelerated inexact Newton method.
In \citet{sleijpen:04} we did this for simple eigenvalue problems using the Jacobi--Davidson method, which is in fact an accelerated Newton method.
The method we present here will make this feasible for the 3D Navier--Stokes equations, even at high Reynolds numbers and close to bifurcation points, for which we have so far failed to find other scalable parallel solution methods.

To achieve this, fully aggregated approaches are important.
In this category, multigrid and multilevel ILU methods for systems of PDEs exist (see \citet{trottenberg:00, wathen:15} and references therein).
The former is attractive, but for those methods smoothers may fail due to a loss of diagonal dominance for higher Reynolds numbers, for which a common solution is to use time integration \citep{feldman:10}.
The latter comprise AMG algorithms and multilevel methods like MRILU \citep{botta:99} and the methods available in ILUPACK \citep{bollhoefer:06}.
ILUPACK has been successful in many fields since it uses a bound to preclude very unstable factorizations.
However, this method does not show grid independence for Navier--Stokes problems and is difficult to parallelize \citep{aliaga:08}.
It works well for large problems, but not yet beyond a single shared memory system.

In this paper, we present a novel multilevel preconditioning method which is specially designed for the 3D Navier--Stokes equations.
In \secref{twolevel}, we first describe the two-level ILU preconditioner as introduced in \citet{wubs:11} and \citet{thies:11b}.
After this, we generalize the two-level method to a multilevel method in \secref{multilevel}.
To make this method work for the 3D Navier--Stokes equations, we introduce a skew partitioning method in \secref{skew}.
In \secref{smilu-results} we discuss the scalability and general performance of the method, and compare it to a popular physics based method, after which we summarize the paper in \secref{smilu-conc}.

\section{The two-level ILU preconditioner}\label{sec:twolevel}
In \citet{wubs:11} a robust parallel two-level method was developed for solving
\begin{align*}
  A\bvec{x}=\bvec{b},
\end{align*}
with $A\in\R^{n\times n}$ for a class of matrices known as $\mathcal{F}$-matrices.
An $\mathcal{F}$-matrix is a matrix of the form
\begin{align*}
  A =
  \begin{pmatrix}
    K & B\\
    B^T & 0
  \end{pmatrix},
\end{align*}
with $K$ symmetric positive definite and $B$ such that it has at most two entries per row and the entries in each row sum up to 0, as is the case for our gradient matrix $G$ \citep{tuma:02,niet:08}.
It has been shown that the two-level preconditioner leads to grid-independent convergence for the Stokes equations on an Arakawa C-grid, and that the method is robust even for the Navier--Stokes equations, which strictly speaking do not yield $\mathcal{F}$-matrices because $K$ becomes nonsymmetric and possibly indefinite.

Rather than reviewing the method and theory in detail, we only briefly present it here.
For details, see \citet{wubs:11} and \citet{thies:11b}.

To simplify the discussion, we focus on the special case of the 3D incompressible Navier--Stokes equations in a cube-shaped domain, discretized on an Arakawa C-grid (see \figref{cgrid}).
We refer to the velocity variables, which are located on the cell faces as $V$-nodes, and to the pressure, which is located in the cell center, as $P$-node.
The variables are numbered cell-by-cell, i.e.\ the first three variables are the $u/v/w$-velocity at the north/east/top face of the cell in the bottom south west corner of the domain, and variable four is the $P$-node in its center.
Appropriate boundary conditions (e.g.\ Dirichlet conditions) are easily implemented in this situation.
We remark that the algorithm (and software) can handle more general situations like rectangular domains, interior boundary cells, etc., and could be generalized to arbitrary domain shapes and partitionings.

First we describe the initialization phase of the preconditioner, which is the necessary preprocessing that has to be done only once for a series of linear systems with matrices with the same sparsity pattern.
Then we describe the factorization phase, which has to be done separately for every matrix.
Finally we describe the solution phase, which has to be performed for every application of the preconditioner.
\subsection{Initialization phase}
\begin{figure}[!ht]
  \begin{center}
    \includegraphics[width=.4\textwidth]{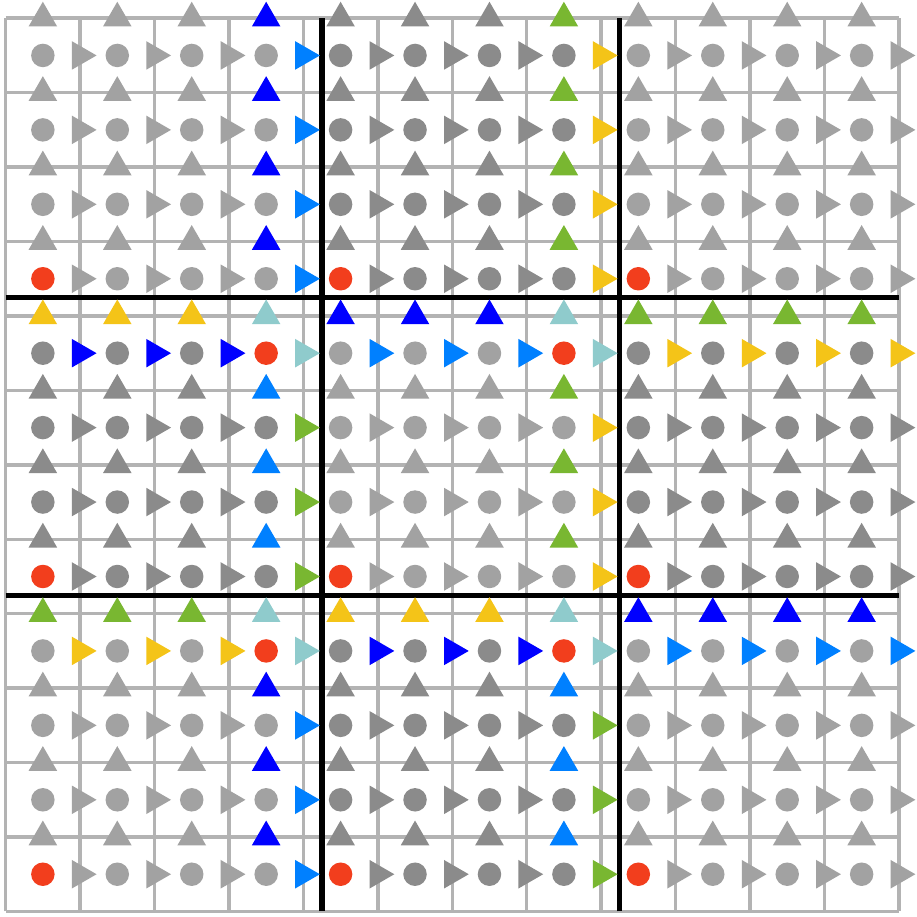}
  \end{center}
  \caption{Cartesian partitioning of a $12 \times 12$ C-grid discretization of the Navier--Stokes equations into 9 subdomains of size $s_x \times s_y$ with $s_x=s_y=4$.
    The interiors are shown in gray.
    Velocities of the same (non-gray) color together form a separator group (but only if they point in the same direction and are in neighboring grid cells).
    The red circles are pressure nodes that are retained.
    \label{fig:cartpart}}
\end{figure}

First the system is partitioned into $N$ non-overlapping subdomains $\Omega_\alpha$, with $\alpha=1,\ldots,N$.
These subdomains may be distributed over different processes, which allows for parallelization of the computation.
The partitioning may be done based on the graph of the matrix, as implemented for instance in METIS \citep{karypis:98}, or by a geometric approach, e.g.\ by dividing the domain into rectangular subdomains.
An example of a Cartesian partitioning of a square domain is shown in \figref{cartpart}.
The non-overlapping subdomains are denoted by the black lines.

Based on the non-overlapping decomposition, we can define a global minimal set of separator nodes $\Gamma$ such that it 
decouples the linear systems associated with the remaining variables in any two subdomains $\Omega_{\alpha, \beta}$. 
This $\Gamma$ can straightforwardly be defined based on either geometric information or the graph of the matrix.
Let $I_\alpha=\Omega_\alpha\backslash\Gamma$ denote the set of interior variables of subdomain $\Omega_\alpha$, then we 
have by construction
\begin{align}
\forall {\alpha\neq\beta, i\in I_\alpha, j\in I_\beta}: A_{i,j}=0, 
\end{align}
Furthermore, we define the set of separator variables associated with subdomain $\Omega_\alpha$ as 
$S_\alpha = \{j\in \Gamma: \exists i \in I_\alpha: A_{i,j}\neq 0\}$.
Note that the interior variables of two subdomains can be eliminated independently on a parallel computer, and that the 
union of the sets $I_\alpha$ and $S_\alpha$ defines an overlapping partitioning. Although we use the non-overlapping 
domain decomposition to define our incomplete factorization, we formally introduce the overlapping subdomains $\overline{\Omega}_\alpha=I_\alpha \cup S_\alpha$.
In the remainder of the paper, we will refer to the variables in the sets $I_\alpha$ and $S_\alpha, \alpha=1,\dots,N$ as \emph{interior nodes}, and \emph{separator nodes}, respectively.
One \emph{separator} is defined as a set of separator nodes that are coupled with the same set of subdomains (geometrically, separators comprise faces, edges and corners).
One \emph{separator group} is defined as the variables on the same separator that have the same variable type ($u,v,w$ or $p$).
Note that because separators are defined by the set of subdomains they are coupled to, separator nodes can never be in multiple separator groups at the same time.
In \figref{cartpart} the interior nodes are shown in gray and the different separator groups are denoted by different colors.

We can now reorder the matrix $A$ such that the interiors ($I$) and separators ($S$) are grouped.
This gives us the system
\begin{align*}
  \begin{pmatrix}
    A_{II} & A_{IS}\\
    A_{SI} & A_{SS}
  \end{pmatrix}
  \begin{pmatrix}
    \bvec{x}_I\\
    \bvec{x}_S
  \end{pmatrix}
  =
  \begin{pmatrix}
    \bvec{b}_I\\
    \bvec{b}_S
  \end{pmatrix},
\end{align*}
where $A_{II}$ consists of independent diagonal blocks.
This submatrix is invertible because on each subdomain we deal with a discretized Navier--Stokes problem on an Arakawa C-grid which is known to be well-posed if the normal and tangential velocities are specified on the boundaries and the level of the pressure is fixed by setting it in one grid point inside the domain.
The velocities are indeed specified on the surrounding boundaries and separators of each subdomain and one pressure is fixed in each subdomain.
Since $A_{II}$ consists of independent diagonal blocks, applying $A_{II}^{-1}$ is easy and trivial to parallelize.
For this reason, we aim to solve
\begin{align*}
  S\bvec{x}_S&=\bvec{b}_S-A_{SI}A_{II}^{-1}\bvec{b}_I,\\
  \bvec{x}_I&=A_{II}^{-1}\bvec{b}_I-A_{II}^{-1}A_{IS}\bvec{x}_S,
\end{align*}
where $S$ is the Schur complement given by $S=A_{SS} - A_{SI}A_{II}^{-1}A_{IS}$.

Whether a variable is coupled to a different subdomain could be determined directly from the graph of the matrix, as discussed earlier.
In our implementation, however, we determine this geometrically by defining the overlapping subdomains during the partitioning step, and checking what overlapping subdomains a node of the non-overlapping subdomain belongs to.

There are three types of separators: faces (in 3D), edges and corners.
For the Navier--Stokes problem on a C-grid, these separators only consist of $V$-nodes.
The $P$-nodes are only connected to $V$-nodes that belong to the same overlapping subdomain, so these should never lie on a separator.
We arbitrarily choose one $P$-node in every interior which we also define to be its own separator group to make sure the Schur complement stays nonsingular.

We want to eliminate velocities on a separator in such a way that the velocities that remain on a separator face provide an approximation of the collective flux through that face.
It is therefore important that the variables are correctly scaled before the factorization in a way that they represent physical fluxes.
In the matrix in \eqref{ns_spp} this gives a $G$-part with entries of constant magnitude.
In this case we can define a Householder transformation which exactly decouples all but one $V$-node from the $P$-nodes \citep{wubs:11}.
This transformation is of the form
\begin{align}
  H_j = I-2\frac{\bvec{v}_j \bvec{v}_j^T}{\bvec{v}_j^T\bvec{v}_j}, \label{eq:householder1}
\end{align}
for some separator group $g_j$ with
\begin{align}
  \bvec{v}_j^{(k)}=
  \begin{cases}
    \bvec{e}_j^{(k)} + \Vert \bvec{e}_j \Vert_2 & \textrm{if node $k$ is the first node of separator group $g_j$}\\
    \bvec{e}_j^{(k)} & \textrm{otherwise}
  \end{cases}\label{eq:householder2}
\end{align}
and
\begin{align*}
  \bvec{e}_j^{(k)}=
  \begin{cases}
    1 & \textrm{if node $k$ is in separator group $g_j$}\\
    0 & \textrm{if node $k$ is not in separator group $g_j$}
  \end{cases}
\end{align*}
for all $k=1,\ldots,n$.
We call the one $V$-node that is still coupled to the $P$-nodes a $V_\Sigma$ node.
Since the Householder transformation can be applied independently for every separator group, we may collect all these transformations in one single transformation matrix $H$.
\subsection{Factorization phase}
For every overlapping subdomain $\overline{\Omega}_i$, $i=1,\ldots,N$, where $N$ is the total number of overlapping subdomains, we can define a local matrix
\begin{align*}
  A^{(i)} =
  \begin{pmatrix}
    A^{(i)}_{II} & A^{(i)}_{IS}\\
    A^{(i)}_{SI} & A^{(i)}_{SS}
  \end{pmatrix}.
\end{align*}
After computing the factorization $A^{(i)}_{II}=L_{II}^{(i)}U_{II}^{(i)}$, the local contribution to the Schur complement is given by
\begin{align*}
  S_i=A_{SS}^{(i)}-A_{SI}^{(i)}(L_{II}^{(i)}U_{II}^{(i)})^{-1}A_{IS}^{(i)},
\end{align*}
and the global Schur complement is given by
\begin{align*}
  S=\sum_{\Omega_i} A_{SS}^{(i)}-\sum_{\overline{\Omega}_i} A_{SI}^{(i)}(L_{II}^{(i)}U_{II}^{(i)})^{-1}A_{IS}^{(i)}.
\end{align*}
Here we take the sum of the $A_{SS}^{(i)}$ contributions over the non-overlapping subdomains to make sure that contributions from separators are not counted multiple times.

We now apply the Householder transformation
\begin{align}
  \begin{split}
    S_T = HSH^T &= H\left(\sum_{\Omega_i} A_{SS}^{(i)}-\sum_{\overline{\Omega}_i} A_{SI}^{(i)}(L_{II}^{(i)}U_{II}^{(i)})^{-1}A_{IS}^{(i)}\right)H^T\\
    &=\sum_{\Omega_i} H_iA_{SS}^{(i)}H_i^T-\sum_{\overline{\Omega}_i} H_iA_{SI}^{(i)}(L_{II}^{(i)}U_{II}^{(i)})^{-1}A_{IS}^{(i)}H_i^T,
  \end{split}
      \label{eq:schurtrans}
\end{align}
which can be done separately for every separator group or subdomain, or on the entire Schur complement.
We choose to apply the transformation separately for every subdomain, since $H$ is very sparse, and sparse matrix-matrix products are very expensive and memory consuming.

We now drop all connections between $V_\Sigma$ and non-$V_\Sigma$ nodes, and between non-$V_\Sigma$ nodes and non-$V_\Sigma$ nodes in different separator groups.
The dropping that is applied here is what makes our factorization inexact.
Not applying the dropping gives us a factorization that can still be partially parallelized, but is also exact.

Our transformed Schur complement is now reduced to a block-diagonal matrix with blocks of non-$V_\Sigma$ nodes for every separator, and one block for all $V_\Sigma$ nodes, which we call $S_{\Sigma\Sigma}$.
Separate factorizations can again be made for all these diagonal blocks, which can again be done in parallel.
For the non-$V_\Sigma$ blocks, a sequential LAPACK solver can be used, and for $S_{\Sigma\Sigma}$ we can employ a (distributed) sparse direct solver.
We denote the factorization that is computed by these solvers as $L_SU_S$.

The full factorization obtained by the method is given by
\begin{align*}
  A \approx
  \begin{pmatrix}
    L_{II} & 0\\
    A_{SI}U_{II}^{-1} & H^TL_S
  \end{pmatrix}
  \begin{pmatrix}
    U_{II} & L_{II}^{-1}A_{IS}\\
    0 & U_SH
  \end{pmatrix}.
\end{align*}

\subsection{Solution phase}
After the preconditioner has been computed, it can be applied in each step of a preconditioned Krylov subspace method, for which we use the well-known GMRES method \citep{saad:86}.
Communication has to occur in the application of $A_{IS}$ and $A_{SI}$, and when solving the distributed $V_\Sigma$ system.
The latter was adressed by using a parallel sparse direct solver in \citet{thies:11b}, but in the next section we propose a different road that does not rely on the availability of such a solver.

\section{The multilevel ILU preconditioner}\label{sec:multilevel}
The main issue with the above two-level ILU factorization that prevents us from scaling up to very large problem sizes in three space dimensions is that, as the problem size increases and the subdomain size remains constant, the size of $S_{\Sigma\Sigma}$ increases steadily, and any (serial or parallel) sparse direct solver eventually limits the feasible problem sizes.
Increasing the subdomain size, on the other hand, leads to more iterations and therefore more synchronization points.

One of the strong points, on the other hand, is the fact that it preserves the structure of the original problem in the sense that, when applied to an $\mathcal{F}$-matrix, it produces a strongly reduced matrix $S_{\Sigma\Sigma}$ which is again an $\mathcal{F}$-matrix.
It is therefore intuitive to try to apply the scheme recursively and avoid the problem of the growing sparse matrix that has to be factorized.
From the structure preserving properties of the algorithm, it is expected that such a recursive scheme again leads to grid-independent convergence if the number of recursions is kept constant as the grid size is increased.

From now on we refer to the number of recursions, or the number of times a reduced matrix $S_{\Sigma\Sigma}$ is computed, as the number of levels.
Note that this means that the method, which was previously referred to the two-level method is in fact the first level of the multilevel method.
Applying a direct solver to $S_T$ from \eqref{schurtrans} is level zero. In this case the preconditioner is a direct solver and GMRES will converge in 1 iteration.

In order to apply the method to the reduced matrix $S_{\Sigma\Sigma}$, we require a coarser partitioning, in which a subdomain consists of multiple subdomains from the original partitioning.
In case we have a regular partitioning like a rectangular partitioning, this may be done by multiplying the \emph{separator length} by a certain \emph{coarsening factor}.
Having a coarsening factor of 2, for instance, means that in 3D the separator length is increased by a factor 2, and the number of subdomains is reduced by a factor 8.

As stated in the previous section, we require the velocity variables to be correctly scaled to be able to apply the Householder transformation.
However, the $V_\Sigma$-variables from the previous level that lie on one separator have a different number of variables in their separator groups.
In case of a regular partitioning, an edge separator, for instance, consists of $V_\Sigma$-nodes from two edges and one corner from the previous level
This leads to a different scaling of the $V_\Sigma$-nodes and thus to non-constant entries in the $G$-part of $S_{\Sigma\Sigma}$.
Instead of re-scaling the $S_{\Sigma\Sigma}$ matrix on every level, we instead use a test vector \bvec{t}.
The test vector is defined initially as a constant vector of ones, and is multiplied by the Householder transformation $H$ at each consecutive level.
The Householder transformation is as defined in \eqref{householder1} and \eqref{householder2}, but now with
\begin{align*}
  \bvec{e}_j^{(k)}=
  \begin{cases}
    \bvec{t}^{(k)} & \textrm{if node $k$ is in separator group $g_j$}\\
    0 & \textrm{if node $k$ is not in separator group $g_j$}
  \end{cases}
\end{align*}
for all $k=1,\ldots,n$.
After applying the Householder transformation, we can again apply dropping to remove connections between $V_\Sigma$ and non-$V_\Sigma$ nodes and between non-$V_\Sigma$ nodes and non-$V_\Sigma$ nodes in different separator groups.
When the matrix $S_{\Sigma\Sigma}$ is sufficiently small, a direct solver is applied to factorize it.

Instead of just having one separator group per variable per separator, we may also choose to have multiple separator groups, meaning that instead of retaining only one $V_\Sigma$ node per variable per separator we retain multiple $V_\Sigma$ nodes.
This in turn means that less dropping is applied, and therefore the factorization is more accurate.
Retaining all nodes in this way, possibly only after reaching a certain level, gives us an exact factorization, which, in terms of iterations for the outer iterative solver, should give the same results as using any other direct solver at that level.
A visual representation of this process is given in \figref{retain}.

\definecolor{sepcol1}{RGB}{121, 183, 49}
\definecolor{sepcol2}{RGB}{0, 128, 255}
\definecolor{sepcol3}{RGB}{0, 0, 255}
\definecolor{sepcol4}{RGB}{244, 196, 24}
\definecolor{sepcol5}{RGB}{143, 203, 204}
\definecolor{sepcol6}{RGB}{242, 62, 29}

\usetikzlibrary{shapes}
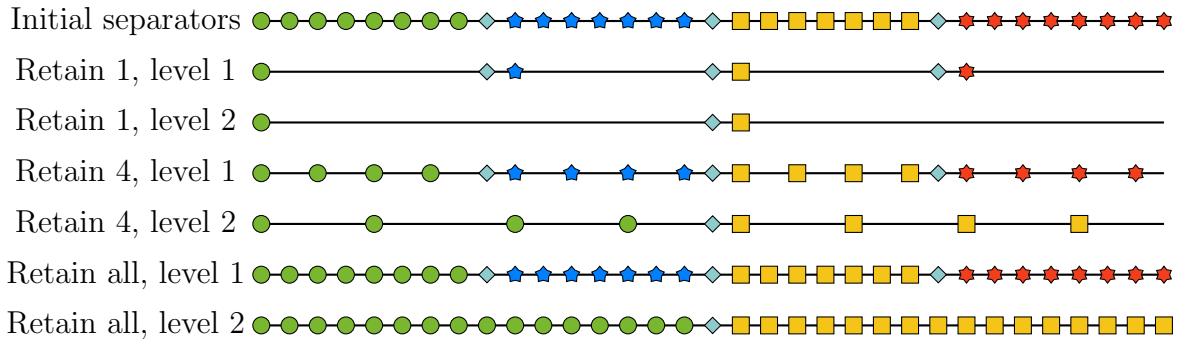
\begin{figure}[ht]
  \centering
  \tikzset{sep1/.style={fill=sepcol1, circle, minimum size=0.15cm, inner sep=0pt, scale=1.5}}
  \tikzset{sep2/.style={fill=sepcol2, star, minimum size=0.15cm, inner sep=0pt, scale=1.5}}
  \tikzset{sep3/.style={fill=sepcol3, rectangle, minimum size=0.15cm, inner sep=0pt, scale=1.5}}
  \tikzset{sep4/.style={fill=sepcol4, rectangle, minimum size=0.15cm, inner sep=0pt, scale=1.5}}
  \tikzset{sep5/.style={fill=sepcol5, diamond, minimum size=0.15cm, inner sep=0pt, scale=1.5}}
  \tikzset{sep6/.style={fill=sepcol6, star, star points=6, minimum size=0.15cm, inner sep=0pt, scale=1.5}}
  \begin{tikzpicture}[scale=1.5, yscale=0.45]
    \draw[thick] (0,6) -- (8,6);
    \foreach \x in {0,...,7} {
      \node[draw, sep1] at (0.25*\x,6) {};
    }
    \node[draw, sep5] at (2,6) {};
    \foreach \x in {1,...,7} {
      \node[draw, sep2] at (2+0.25*\x,6) {};
    }
    \node[draw, sep5] at (4,6) {};
    \foreach \x in {1,...,7} {
      \node[draw, sep4] at (4+0.25*\x,6) {};
    }
    \node[draw, sep5] at (6,6) {};
    \foreach \x in {1,...,8} {
      \node[draw, sep6] at (6+0.25*\x,6) {};
    }
    \node at (-1.2,6){Initial separators};

    \draw[thick] (0,5) -- (8,5);
    \foreach \x in {0,...,0} {
      \node[draw, sep1] at (0.25*\x,5) {};
    }
    \node[draw, sep5] at (2,5) {};
    \foreach \x in {1,...,1} {
      \node[draw, sep2] at (2+0.25*\x,5) {};
    }
    \node[draw, sep5] at (4,5) {};
    \foreach \x in {1,...,1} {
      \node[draw, sep4] at (4+0.25*\x,5) {};
    }
    \node[draw, sep5] at (6,5) {};
    \foreach \x in {1,...,1} {
      \node[draw, sep6] at (6+0.25*\x,5) {};
    }
    \node at (-1.2,5){Retain 1, level 1};

    \draw[thick] (0,4) -- (8,4);
    \foreach \x in {0,...,0} {
      \node[draw, sep1] at (2.25*\x,4) {};
    }
    \node[draw, sep5] at (4,4) {};
    \foreach \x in {0,...,0} {
      \node[draw, sep4] at (4.25+2*\x,4) {};
    }
    \node at (-1.2,4){Retain 1, level 2};

    \draw[thick] (0,3) -- (8,3);
    \foreach \x in {0,...,3} {
      \node[draw, sep1] at (0.5*\x,3) {};
    }
    \node[draw, sep5] at (2,3) {};
    \foreach \x in {0,...,3} {
      \node[draw, sep2] at (2.25+0.5*\x,3) {};
    }
    \node[draw, sep5] at (4,3) {};
    \foreach \x in {0,...,3} {
      \node[draw, sep4] at (4.25+0.5*\x,3) {};
    }
    \node[draw, sep5] at (6,3) {};
    \foreach \x in {0,...,3} {
      \node[draw, sep6] at (6.25+0.5*\x,3) {};
    }
    \node at (-1.2,3){Retain 4, level 1};

    \draw[thick] (0,2) -- (8,2);
    \foreach \x in {0,...,1} {
      \node[draw, sep1] at (\x,2) {};
    }
    \foreach \x in {0,...,1} {
      \node[draw, sep1] at (2.25+\x,2) {};
    }
    \node[draw, sep5] at (4,2) {};
    \foreach \x in {0,...,3} {
      \node[draw, sep4] at (4.25+\x,2) {};
    }
    \node at (-1.2,2){Retain 4, level 2};

    \draw[thick] (0,1) -- (8,1);
    \foreach \x in {0,...,7} {
      \node[draw, sep1] at (0.25*\x,1) {};
    }
    \node[draw, sep5] at (2,1) {};
    \foreach \x in {1,...,7} {
      \node[draw, sep2] at (2+0.25*\x,1) {};
    }
    \node[draw, sep5] at (4,1) {};
    \foreach \x in {1,...,7} {
      \node[draw, sep4] at (4+0.25*\x,1) {};
    }
    \node[draw, sep5] at (6,1) {};
    \foreach \x in {1,...,8} {
      \node[draw, sep6] at (6+0.25*\x,1) {};
    }
    \node at (-1.2,1){Retain all, level 1};

    \draw[thick] (0,0) -- (8,0);
    \foreach \x in {0,...,15} {
      \node[draw, sep1] at (0.25*\x,0) {};
    }
    \node[draw, sep5] at (4,0) {};
    \foreach \x in {1,...,16} {
      \node[draw, sep4] at (4+0.25*\x,0) {};
    }
    \node at (-1.2,0){Retain all, level 2};
  \end{tikzpicture}
  \caption{One-dimensional representation of the process of retaining multiple nodes per separator.
  Each separator group has its own color and shape.}
\label{fig:retain}
\end{figure}
\section{Skew partitioning in 2D and 3D}\label{sec:skew}
Looking at \figref{cartpart}, we see that there are pressures that are located at the corners of the subdomains that are surrounded by velocity separators.
This means that if we add these pressures to the interior, as was suggested above, we get a singular interior block.
We call these pressure nodes that are surrounded by velocity separators \emph{isolated pressure nodes}.
For the two-level preconditioner in 2D, it was possible to solve this problem by adding these pressures to the Schur complement as single-element separator groups.
This unfortunately does not work for the multilevel method, since in this case velocity nodes around the isolated pressure nodes get eliminated.
It also does not work in 3D because there the isolated pressure nodes also exist inside of the edge separators, where they form tubes of isolated nodes.

A possible way to solve this for the multilevel case and in 3D is to also turn all velocity nodes around the isolated pressure nodes into separate separator groups.
This means that they are all treated as $V_\Sigma$ nodes and are never eliminated until they are in the interior of the domain at a later level.
This, however, has the downside that a lot more nodes have to be retained at every level, resulting in much larger $S_{\Sigma\Sigma}$ matrices at every level.
Furthermore, a lot of bookkeeping is required to keep track of all the extra separator groups, and their elimination leads to long-range connections and thus a more irregular matrix structure.

\begin{figure}[bpt]
  \centering
  \subtop[Initial partitioning\label{fig:skewparta}]{
    \includegraphics[width=0.35\textwidth]{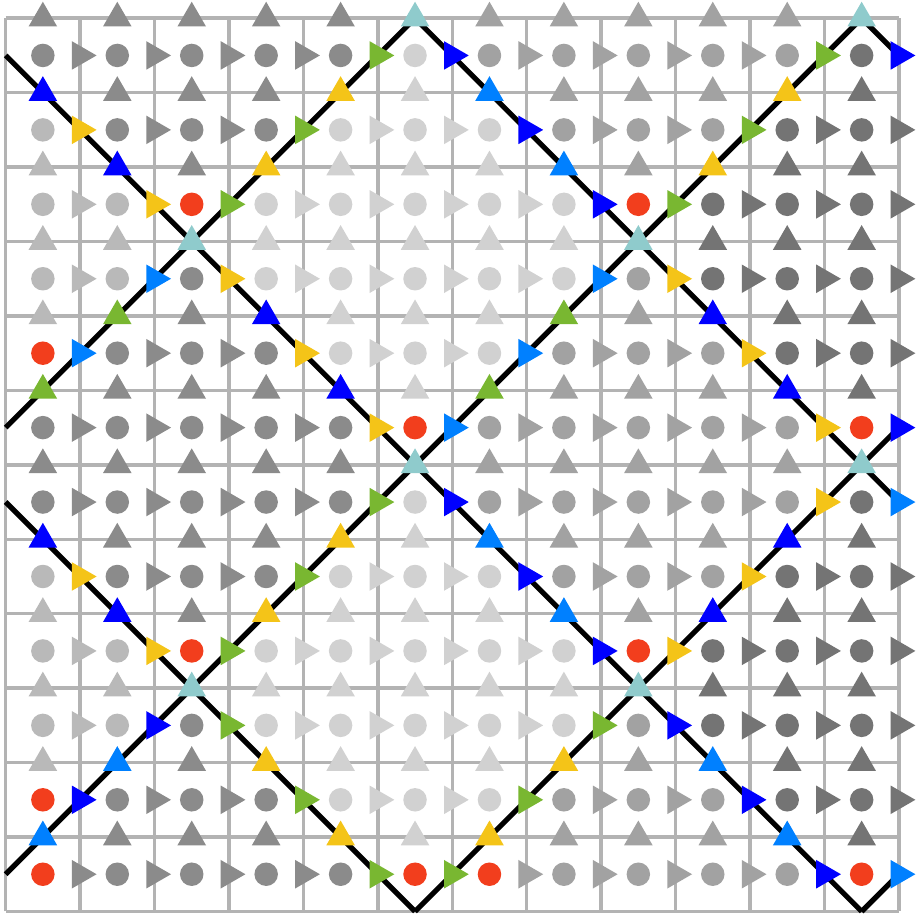}}
  \subtop[After elimination of interiors only the separator groups remain]{
    \includegraphics[width=0.35\textwidth]{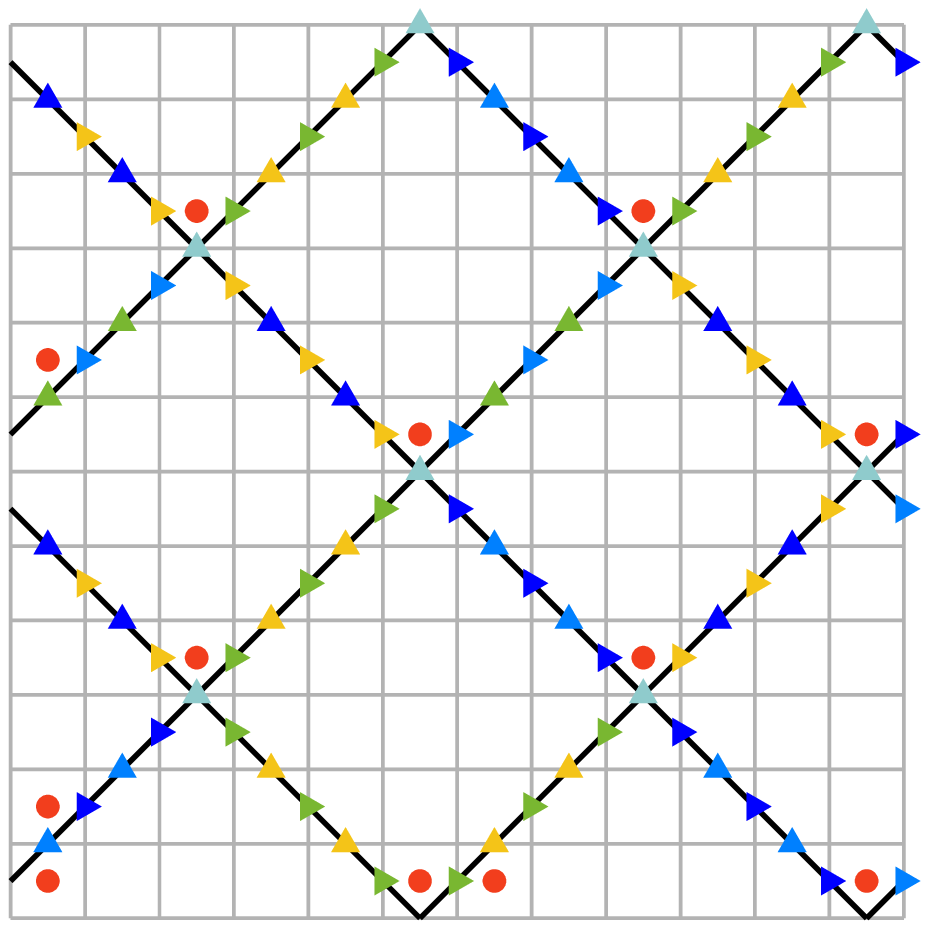}}\\
  \subtop[The Householder transformations leave one $V_\Sigma$-node per separator group]{
    \includegraphics[width=0.35\textwidth]{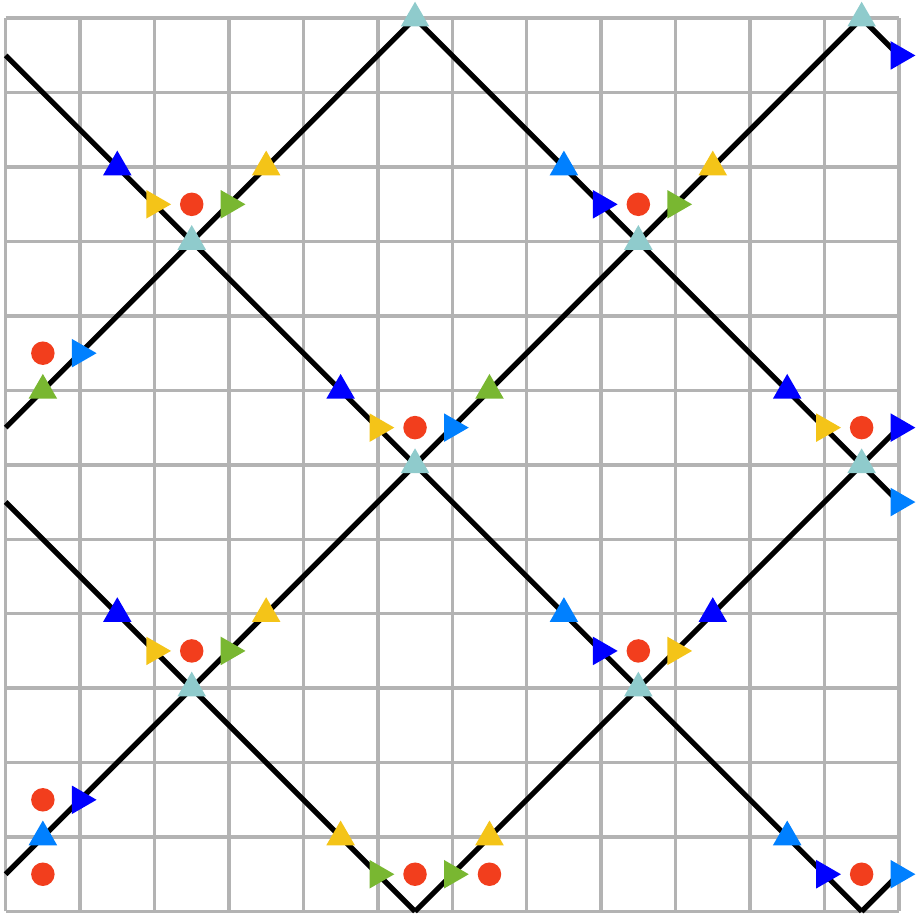}}
  \subtop[Partitioning on next level]{
    \includegraphics[width=0.35\textwidth]{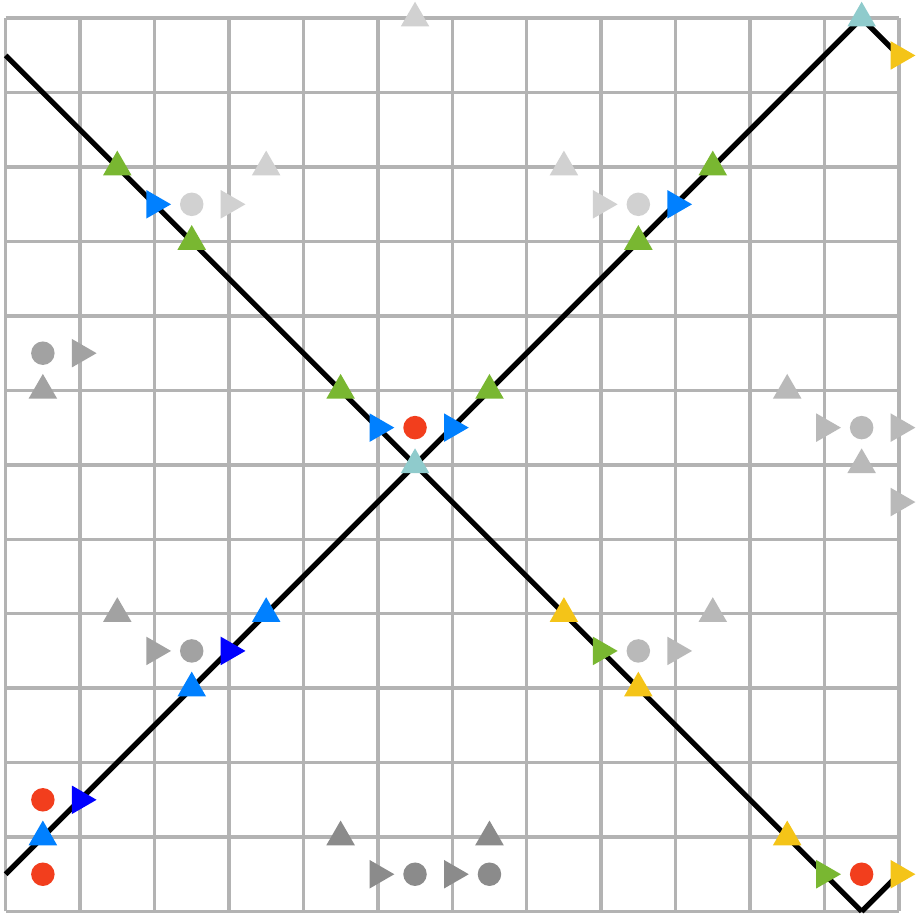}}
  \caption{Skew partitioning of a $12 \times 12$ C-grid discretization of the Navier--Stokes equations into 12 subdomains.
    The interiors are shown in gray.
    Velocities of the same (non-gray) color together form a separator group (but only if they point in the same direction and are in neighboring grid cells).
    The red circles are pressure nodes that are retained.
    This example uses a coarsening factor of 2, i.e.\ the separators on the next level have twice the length of those on the previous level.}
\label{fig:skewpart}
\end{figure}

In 2D, we can resolve these problems by rotating the Cartesian partitioning by 45 degrees.
The result is shown in \figref{skewparta}.
It is easy to see that in this case, no pressure nodes are surrounded by only velocity separators.
We call this partitioning method \emph{skew partitioning}.
In \figref{skewpart}, we also show the workings of the multilevel method, with all the steps being displayed in the different subfigures.

For the skew partitioning to work with our multilevel method, we have two requirements on the shape of the subdomains: (1) it should be \emph{space-filling}, meaning that it can be used to fill the entire domain and (2) it should be \emph{self-similar}, meaning that we can construct a larger subdomain of the same shape from multiple smaller subdomains.
It is easy to see that these two properties hold for the 2D skew partitioning.

The most basic idea for generalizing the rotated square shape to a 3D setting is to use octahedral subdomains.
Partitioning with this design turned out to be unsuccessful, but it is still briefly discussed here since it led to some new insights.
Since regular octahedra (the dual to cubes, having their vertices at the centers of the cube faces) by themselves are not space-filling, the octahedra can be slightly distorted to make them fit within a single cubic repeat unit.
The resulting subdomains are space-filling, but only by using three mutually orthogonal subdomain types.
The fact that it requires the use of three types of subdomains increases the programming efforts significantly since it introduces a lot of edge cases that should be considered, e.g.\ for subdomains located at the boundary of the domain.

The major problem with the octahedral subdomains, however, is that they are not self-similar, meaning that we cannot construct a larger octahedral subdomain from multiple smaller octahedral subdomains.
However, self-similarity can be achieved by splitting the octahedra into four smaller tetrahedra, of which six different types are required to fill 3D space.
This introduces additional separator planes that are similar to the 2D skew case and hence it increases the risk of isolating a pressure node when such planes intersect.
Especially planar intersections which are parallel to any of the Cartesian axes have a high risk of producing isolated pressure nodes.

\begin{figure}[!ht]
  \begin{center}
    \subtop[Shape of the parallelepiped inside of two cubes that consist of $m\times m\times m$ grid 
    cells\label{fig:ppda}]{
      \includegraphics[width=0.3\textwidth]{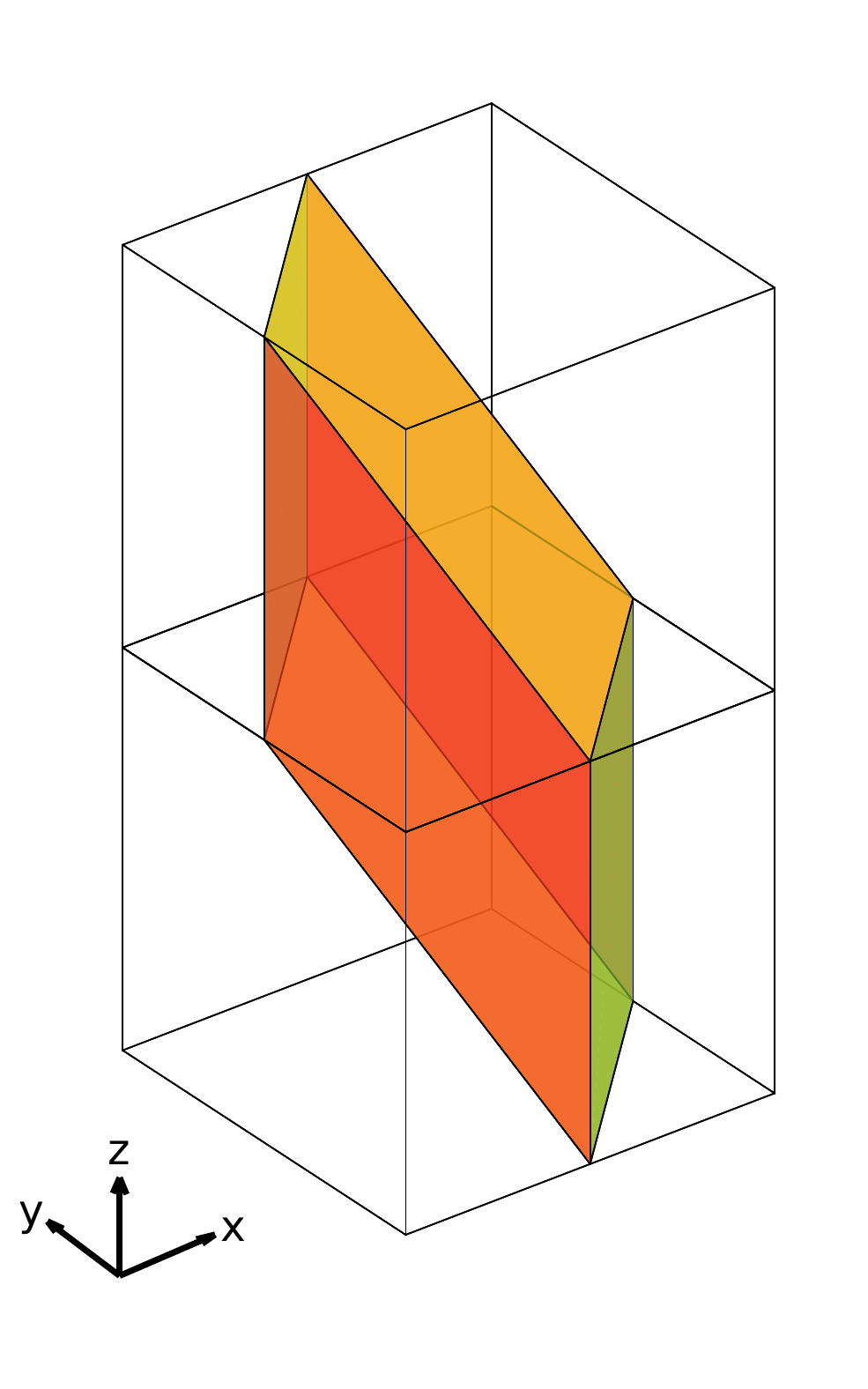}}
    \hspace{1em}
    \subtop[Schematic view of the position of the separator nodes from the same point of view as the figure on the left The $u$-nodes are depicted as red faces, the $v$-nodes are depicted in green and the $w$-nodes in yellow.\label{fig:ppdb}]{
      \includegraphics[width=0.3\textwidth]{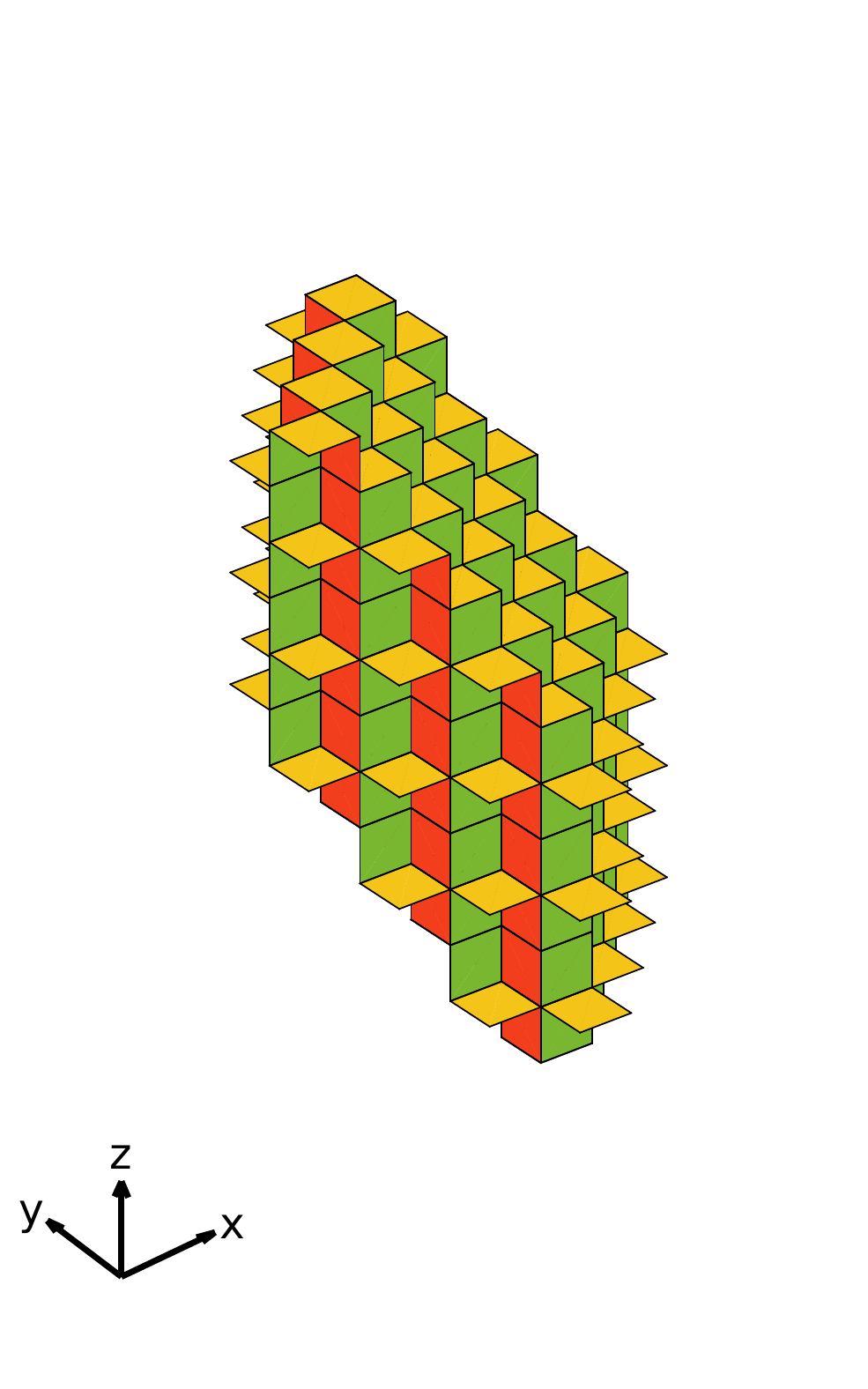}}
  \end{center}
  \caption{Parallelepiped shape for partitioning the domain.} \label{fig:ppd}
\end{figure}

We did indeed not manage to find any self-similar tiling using tetrahedra that would not give rise to any isolated pressure nodes.
Moreover, we would like to have a single subdomain shape that we can use instead of six, since this would greatly simplify the implementation.
A lesson we learned is that isolated pressure nodes always seem to arise when having faces that are aligned with the grid.
Therefore, we looked into a rotated parallelepiped shape that does not have any faces that are aligned with the grid \citep{klok:17}.
This shape is shown in \figref{ppda}, where the cubes represent a set of $s_x\times s_x\times s_x$ grid cells.
A welcome property of this domain is that its central cross section resembles the proposed 2D skew partitioning method.

A schematic view of the position of the separator nodes is shown in \figref{ppdb}.
One important detail to note is that on the side that is facing towards us, only half of the $w$-nodes are displayed.
This is because the $w$-nodes have to be positioned in an alternating fashion on the inside and outside of the non-overlapping subdomain to prevent isolated pressure nodes from appearing.
A consequence of this alternating property is that the $w$-planes have to be divided into two separate separator groups; one for the $w$-nodes that are located inside the non-overlapping subdomain, and one for $w$-nodes that are only present in the overlapping subdomain.

Another advantage of using a skew domain partitioning is that the amount of communication that is required is reduced when compared to a square partitioning.
In \citet{bisseling:04}, it is estimated that for the Laplace problem, the communication is asymptotically reduced by a factor of $\sqrt{2}$ for the 2D diamond shape.
If we instead compare the diamond shape to a rectangular domain with the same number of nodes (having the same number of nodes with a square domain is impossible), we find that communication is reduced by a factor of $\frac{3}{2}$.
In the same manner, we can compare a 3D domain of size $s_x\times s_x \times s_x/2$ to the rotated parallelepiped, and find a factor of $\frac{4}{3}$.
We remark that the truncated octahedron that is used in \citet{bisseling:04} for the 3D domain and has a factor of 1.68 cannot be used for our multilevel method, since truncated octahedra are not self-similar.

\newcommand{\Rey}{\text{Re}}
\section{Numerical results}\label{sec:smilu-results}
A common benchmark in fluid dynamics is the lid-driven cavity problem.
We consider an incompressible Newtonian fluid in a square three-dimensional domain of unit length, with a lid at the top which moves at a constant speed $U$.
The equations are given by \eqref{smilu-ns}.
No-slip boundary conditions are applied at the walls, meaning that they are Dirichlet boundary conditions, and the equations are discretized on an Arakawa C-grid as described before.
We take $n_x=n_y=n_z$ grid cells in every direction.

At first, however, we will only look at the Stokes equations of the form
\begin{align*}
  \Delta \bvec{u} - \nabla p &= \bvec{f},\\
  \nabla \cdot \bvec{u} &= 0,
\end{align*}
where we take $\bvec{f}$ to be random.
This is because our preconditioner is constructed in such a way that memory usage and time cost for both computation and application of the preconditioner should not be influenced by inclusion of the convective term.
After this, we further investigate the behavior of the method on the lid-driven cavity problem for increasing Reynolds numbers, which constitute harder problems.
Therefore we expect an increase in iterations of the iterative solver, but otherwise the same behavior.

For obtaining the exact memory usage, we developed a custom library which overrides all memory allocation routines when linking against it.
The library contains a hash table in which the amount of memory that is allocated is stored by its memory address.
We keep track of the total amount of memory that is allocated, which is increased on memory allocation, and reduced by the amount that is stored in the hash table when memory is freed.
The reason we developed this library is that existing methods rely on rough estimates of the memory usage of the used data structures, use the data that is available from \texttt{/proc/meminfo}, which is inaccurate, or actually count memory usage in a similar way as we do (e.g.\ valgrind), but have a large performance impact.

We perform every experiment twice: once to determine the memory usage, and once to determine the run time, without linking to the memory usage library.
This means that when reporting timing results, we are not affected by the performance impact of tools to determine memory usage.
The reason that we are still concerned about their performance impact, even though we developed our own library, is that it adds roughly a constant amount of time per process, which impacts scalability results.
The memory usage that we report is the exact difference in memory usage before and after a certain action is performed, e.g.\ before and after the construction of the preconditioner.

For the timing results, we do not include the time it takes to compute the partitioning, because we did not optimize this step.
The partitioning is computed by first constructing one full overlapping subdomain of the prescribed subdomain size sequentially, which is then mapped to the correct position of every overlapping subdomain to determine the interiors and separator groups.
However, since the initial full overlapping subdomain contains all (possibly already eliminated) nodes, every time we increase the number of levels by 1, the amount of time required to compute the partitioning increases by a factor 8, which is the worst possible scenario.
The reason we do not include this in the timing results is that this may be resolved by only computing the partitioning for nodes that are still present in the Schur complement at that level.
This requires a full rewrite of the existing partitioning code, and will be addressed in a future version of the software.
It does, however, not have any impact on the timing results of the remainder of the implementation, since this is completely decoupled.
This means that even though the partitioning method does not scale at all, the preconditioning method itself can be studied reliably without including the timing of the partitioning.

For the implementation of the preconditioner and solver, we use libraries from the Trilinos project \citep{heroux:05}.
The libraries we use are Epetra (vector and matrix data structures), IFPACK (direct solver and preconditioning interfaces) \citep{sala:05}, Amesos (direct solvers) \citep{sala:08} and Belos (iterative solvers) \citep{bavier:12}.
As iterative solver we use GMRES(250) \citep{saad:86}, as parallel sparse direct solver on the coarsest level we use SuperLU\_DIST 6.1 \citep{liu:18}, and as direct solver for the interior blocks we use KLU \citep{davis:10} with the fill-reducing ordering from \citet{niet:08}.
The implementation of our preconditioner can be found on GitHub~\footnote{\url{https://github.com/nlesc-smcm/hymls}}.

The benchmark is performed on the SuperMUC-NG cluster at LRZ, Munich~\footnote{\url{https://www.top500.org/system/179566/}}.
The nodes contain two Intel Xeon ``Skylake'' Platinum 8174 processors with 24 cores and have 96 GB of memory.
For all experiments, we disable multithreading inside the MPI processes, since the use of OpenMP in Epetra is not compatible with our need for very small objects (e.g.\ sparse matrices living on a very small subdomain).

\subsection{Weak scalability}
First, we look at results obtained when increasing the grid size $n_x$ at the same rate as the number of used cores $n_p$, i.e.\ the problem size on each core is kept constant.
The exact configurations that we use are 1 core for $n_x=16$, 1 core for $n_x=32$, 8 cores for $n_x=64$, 64 cores for $n_x=128$, 512 cores for $n_x=256$ and 4096 cores for $n_x=512$.
The size of the subdomains (the size of the cubes in \figref{ppda}) at the first level is $s_x=8$, while we choose the coarsening factor to be 2, meaning that we increase $s_x$ by a factor of 2 at each level.
We perform experiments where we keep the number of levels constant at $L=2$, and experiments where we increase the number of levels by 1 when doubling the grid size.
For the latter, we look at three cases where we retain a different number of separator nodes starting at level 2: 1, 4, and all.
Since we start retaining more nodes after two levels, results with only two levels ($n_x=16$) will be the same for all configurations.
When no result is shown for $n_x=512$, this means the method ran out of memory, unless otherwise stated.

For the fixed number of levels, we expect the number of iterations of the iterative solver to converge to a constant number as the grid is refined.
For the case where we increase the number of levels as the domain size increases, we expect the number of iterations to only increase mildly, and we expect that retaining more separator nodes starting at level 2 decreases the number of iterations until it again becomes constant as we retain more and more nodes per separator.

\setlength{\figureWidth}{0.5\textwidth}
\setlength{\figureHeight}{0.3\textwidth}
\newcommand\figureScale{1.0}

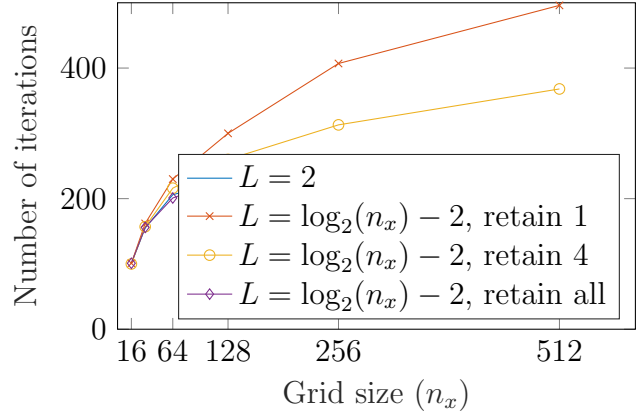
\begin{figure}[!ht]
  \begin{center}
%
%
\definecolor{mycolor1}{rgb}{0.00000,0.44700,0.74100}%
\definecolor{mycolor2}{rgb}{0.85000,0.32500,0.09800}%
\definecolor{mycolor3}{rgb}{0.92900,0.69400,0.12500}%
\definecolor{mycolor4}{rgb}{0.49400,0.18400,0.55600}%
\begin{tikzpicture}[%
scale=\figureScale
]

\begin{axis}[%
width=0.951\figureWidth,
height=\figureHeight,
at={(0\figureWidth,0\figureHeight)},
scale only axis,
xmin=0,
xmax=600,
xtick={ 16,  64, 128, 256, 512},
xlabel style={font=\color{white!15!black}},
xlabel={Grid size ($n_x$)},
ymin=0,
ymax=500,
ylabel style={font=\color{white!15!black}},
ylabel={Number of iterations},
axis background/.style={fill=white},
legend style={at={(0.97,0.03)}, anchor=south east, legend cell align=left, align=left, draw=white!15!black}
]
\addplot [color=mycolor1]
  table[row sep=crcr]{%
16	101\\
32	155\\
64	206\\
128	229\\
256	238\\
};
\addlegendentry{$L = 2$}

\addplot [color=mycolor2, mark=x, mark options={solid, mycolor2}]
  table[row sep=crcr]{%
16	101\\
32	162\\
64	230\\
128	300\\
256	407\\
512	496\\
};
\addlegendentry{$L = \log_2(n_x)-2$, retain 1}

\addplot [color=mycolor3, mark=o, mark options={solid, mycolor3}]
  table[row sep=crcr]{%
16	100\\
32	157\\
64	216\\
128	260\\
256	313\\
512	368\\
};
\addlegendentry{$L = \log_2(n_x)-2$, retain 4}

\addplot [color=mycolor4, mark=diamond, mark options={solid, mycolor4}]
  table[row sep=crcr]{%
16	101\\
32	156\\
64	201\\
128	231\\
256	240\\
};
\addlegendentry{$L = \log_2(n_x)-2$, retain all}

\end{axis}
\end{tikzpicture}%
  \end{center}
  \caption{Number of iterations of GMRES for the Stokes problem on a grid of size $n_x \times n_x \times n_x$, while keeping the number of levels at $L=2$, and when increasing the number of levels by 1 when $n_x$ is doubled.
    A relative residual of $10^{-8}$ was used as convergence tolerance.} \label{fig:smilu-its}
\end{figure}

The results are shown in \figref{smilu-its}.
We see that indeed the number of iterations becomes constant when fixing the number of levels or when retaining all separator nodes.
When increasing the number of levels with the grid size, we see that the number of iterations keeps increasing gradually.
What is interesting is that when we retain 4 instead of 1 separator node per separator group after level 2, the number of iterations that is required decreases drastically, even though this does not have a significant impact on the memory usage as we will see later.

\begin{figure}[!ht]
  \begin{center}
%
%
\definecolor{mycolor1}{rgb}{0.00000,0.44700,0.74100}%
\definecolor{mycolor2}{rgb}{0.85000,0.32500,0.09800}%
\definecolor{mycolor3}{rgb}{0.92900,0.69400,0.12500}%
\definecolor{mycolor4}{rgb}{0.49400,0.18400,0.55600}%
\begin{tikzpicture}[%
scale=\figureScale
]

\begin{axis}[%
width=0.951\figureWidth,
height=\figureHeight,
at={(0\figureWidth,0\figureHeight)},
scale only axis,
xmin=0,
xmax=600,
xtick={ 16,  64, 128, 256, 512},
xlabel style={font=\color{white!15!black}},
xlabel={Grid size ($n_x$)},
ymin=0,
ymax=50,
ylabel style={font=\color{white!15!black}},
ylabel={Time (s)},
axis background/.style={fill=white},
legend style={legend cell align=left, align=left, draw=white!15!black}
]
\addplot [color=mycolor1]
  table[row sep=crcr]{%
16	0.65\\
32	6.42\\
64	9.86\\
128	11.86\\
256	23.92\\
};
\addlegendentry{$L = 2$}

\addplot [color=mycolor2, mark=x, mark options={solid, mycolor2}]
  table[row sep=crcr]{%
16	0.65\\
32	7.05\\
64	9.84\\
128	11.23\\
256	12.59\\
512	16.28\\
};
\addlegendentry{$L = \log_2(n_x)-2$, retain 1}

\addplot [color=mycolor3, mark=o, mark options={solid, mycolor3}]
  table[row sep=crcr]{%
16	0.64\\
32	6.37\\
64	10.05\\
128	11.95\\
256	14.19\\
512	18.14\\
};
\addlegendentry{$L = \log_2(n_x)-2$, retain 4}

\addplot [color=mycolor4, mark=diamond, mark options={solid, mycolor4}]
  table[row sep=crcr]{%
16	0.64\\
32	6.4\\
64	11.06\\
128	44.77\\
256	1415.64\\
};
\addlegendentry{$L = \log_2(n_x)-2$, retain all}

\end{axis}
\end{tikzpicture}%
  \end{center}
  \caption{Time to compute the preconditioner for the Stokes problem on a grid of size $n_x \times n_x \times n_x$, while keeping the number of levels at $L=2$, and when increasing the number of levels by 1 when $n_x$ is doubled.} \label{fig:smilu-time-compute}
\end{figure}
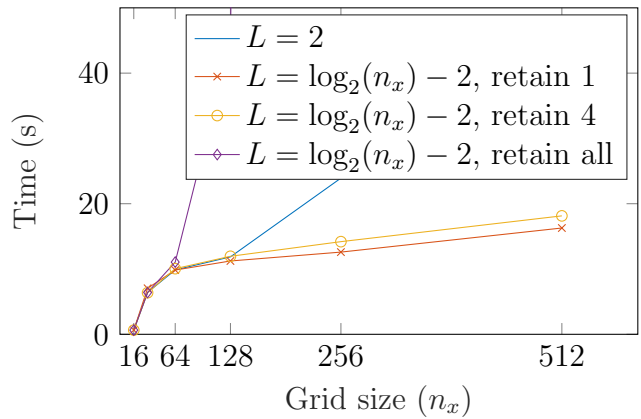

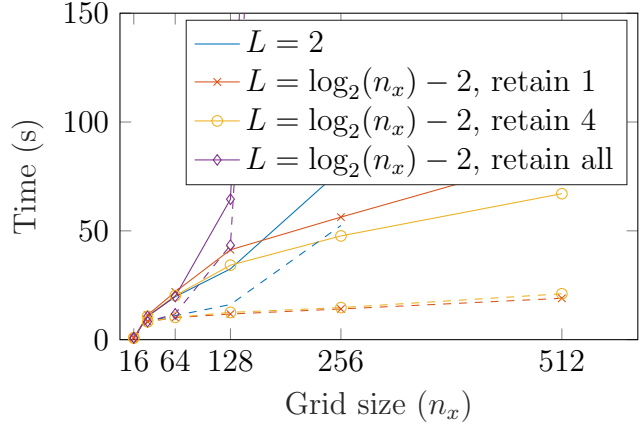
\begin{figure}[!ht]
  \begin{center}
%
%
\definecolor{mycolor1}{rgb}{0.00000,0.44700,0.74100}%
\definecolor{mycolor2}{rgb}{0.85000,0.32500,0.09800}%
\definecolor{mycolor3}{rgb}{0.92900,0.69400,0.12500}%
\definecolor{mycolor4}{rgb}{0.49400,0.18400,0.55600}%
\begin{tikzpicture}[%
scale=\figureScale
]

\begin{axis}[%
width=0.951\figureWidth,
height=\figureHeight,
at={(0\figureWidth,0\figureHeight)},
scale only axis,
xmin=0,
xmax=600,
xtick={ 16,  64, 128, 256, 512},
xlabel style={font=\color{white!15!black}},
xlabel={Grid size ($n_x$)},
ymin=0,
ymax=150,
ylabel style={font=\color{white!15!black}},
ylabel={Time (s)},
axis background/.style={fill=white},
legend style={legend cell align=left, align=left, draw=white!15!black}
]
\addplot [color=mycolor1]
  table[row sep=crcr]{%
16	0.7\\
32	10.87\\
64	19.73\\
128	32.52\\
256	77.96\\
};
\addlegendentry{$L = 2$}

\addplot [color=mycolor2, mark=x, mark options={solid, mycolor2}]
  table[row sep=crcr]{%
16	0.7\\
32	11.36\\
64	21.97\\
128	41.29\\
256	56.31\\
512	85.09\\
};
\addlegendentry{$L = \log_2(n_x)-2$, retain 1}

\addplot [color=mycolor3, mark=o, mark options={solid, mycolor3}]
  table[row sep=crcr]{%
16	0.67\\
32	10.73\\
64	20.22\\
128	34.23\\
256	47.67\\
512	67.07\\
};
\addlegendentry{$L = \log_2(n_x)-2$, retain 4}

\addplot [color=mycolor4, mark=diamond, mark options={solid, mycolor4}]
  table[row sep=crcr]{%
16	0.68\\
32	10.69\\
64	19.87\\
128	64.53\\
256	1112.35\\
};
\addlegendentry{$L = \log_2(n_x)-2$, retain all}

\addplot [color=mycolor1, dashed, forget plot]
  table[row sep=crcr]{%
16	0.95\\
32	8.44\\
64	11.04\\
128	16.15\\
256	52.47\\
};
\addplot [color=mycolor2, dashed, mark=x, mark options={solid, mycolor2}, forget plot]
  table[row sep=crcr]{%
16	0.95\\
32	8.28\\
64	10.11\\
128	11.78\\
256	14.1\\
512	19.01\\
};
\addplot [color=mycolor3, dashed, mark=o, mark options={solid, mycolor3}, forget plot]
  table[row sep=crcr]{%
16	0.93\\
32	8.19\\
64	10.25\\
128	12.53\\
256	14.74\\
512	21.08\\
};
\addplot [color=mycolor4, dashed, mark=diamond, mark options={solid, mycolor4}, forget plot]
  table[row sep=crcr]{%
16	0.93\\
32	8.28\\
64	11.85\\
128	43.4\\
256	890.72\\
};
\end{axis}
\end{tikzpicture}%
  \end{center}
  \caption{Time to solve the linear system with GMRES (lines), and time of 200 applications of the preconditioner (dashed lines).
    This is for the Stokes problem on a grid of size $n_x \times n_x \times n_x$, while keeping the number of levels at $L=2$, and when increasing the number of levels by 1 when $n_x$ is doubled.
    A relative residual of $10^{-8}$ was used as convergence tolerance.}
\label{fig:smilu-time-solve}
\end{figure}

The computational time of both computing the preconditioner, as well as the application of the preconditioner would ideally become constant when keeping the problem size at each core constant while increasing the problem size.
However, in practice this is not possible since increasing the number of cores also increases the required amount of communication.
The results are shown in \figref{smilu-time-compute} and \figref{smilu-time-solve}.

When computing the preconditioner, we see that if we increase the number of levels with the grid size, the computational time rises only slightly.
Ideally, there would be no rise at all, but since an increased the number of computational nodes also means more communication, there will always be an increase in practice.
In our case, the communication mainly happens at the point where contributions of neighboring subdomains have to be added up in the Schur complement.
Since retaining more nodes per level means an increase in the amount of communication, we also see that retaining 4 or all nodes takes more time than retaining only 1 node per separator at every level.
It must be noted that due to system variability, we also found results that were about 10\% better or worse than the results displayed here.
This could be improved by disabling frequency adaptation in the CPUs, but we are aware that at very large core counts the synchronous nature of our algorithm may become a scalability problem.

If we keep the number of levels constant, the computational time required to compute the factorization at the last level increases rapidly, and for larger grid sizes the amount of memory that is required for the factorization is too large for the computation to complete.
We also notice that retaining all nodes after level 2 is much less efficient than just using SuperLU\_DIST at level 2, meaning that our preconditioner performs very poorly as a direct solver.
This is mainly due to the fact that the Schur complement at the last level is quite large and full.
The last level Schur complement for a grid of size $256^3$ has size $20\,961\times 20\,961$, and its factorization is filled with 72\% nonzeros.
Computing the factorization of this matrix and using it in a forward/backward substitution is therefore very expensive.
Using a parallel dense solver instead of SuperLU\_DIST might help to make it more efficient.
Moreover, since we choose $n\sim n_p$, the cost of computing the factorization by a sparse direct solver grows at best with
$\mathcal{O}(n^2/n_p)=\mathcal{O}(n)=\mathcal{O}(n_x^3)$ \cite[]{liu:18}, and we expect that the cost of computing the preconditioner when keeping the number of levels constant, or when retaining all separator nodes, increases with $n_x^3$.

In \figref{smilu-time-solve}, we show both the time required to solve the linear system after computation of the preconditioner, and the time of 200 applications of the preconditioner, which is not influenced by the number of iterations of the iterative solver and does not include the time spent on, for instance, matrix-vector products and orthogonalization.
First of all, we again observe that retaining all separator nodes is a bad idea, since the computational time goes off the chart.
For the case where we only use 2 levels, we also see the unwanted behavior of a time that keeps increasing linearly for the total solution time, and superlinearly for the application of the preconditioner, where we would actually expect $\mathcal{O}(n^{4/3}/n_p)=\mathcal{O}(n_x)$ behavior from the triangular solve \citep{liu:18}.
This may be caused by disabling multithreading support, which results in an increased amount of communication inside of SuperLU\_DIST.

For both cases where we retain 1 and 4 separator nodes after 2 levels, we see that the computational time only slightly increases for larger grid sizes.
We also note that for the case where we retain 4 nodes, the application of the preconditioner is slightly slower, but the total solution time is much smaller due to the lower number of iterations that is required, as can also be seen in \figref{smilu-its}, and the fact that applying the preconditioner is relatively cheap.

\begin{figure}[!ht]
  \begin{center}
%
%
\definecolor{mycolor1}{rgb}{0.00000,0.44700,0.74100}%
\definecolor{mycolor2}{rgb}{0.85000,0.32500,0.09800}%
\definecolor{mycolor3}{rgb}{0.92900,0.69400,0.12500}%
\definecolor{mycolor4}{rgb}{0.49400,0.18400,0.55600}%
\begin{tikzpicture}[%
scale=\figureScale
]

\begin{axis}[%
width=0.951\figureWidth,
height=\figureHeight,
at={(0\figureWidth,0\figureHeight)},
scale only axis,
xmin=0,
xmax=600,
xtick={ 16,  64, 128, 256, 512},
xlabel style={font=\color{white!15!black}},
xlabel={Grid size ($n_x$)},
ymin=0,
ymax=450,
ylabel style={font=\color{white!15!black}},
ylabel={Memory usage (MB)},
axis background/.style={fill=white},
legend style={at={(0.97,0.03)}, anchor=south east, legend cell align=left, align=left, draw=white!15!black}
]
\addplot [color=mycolor1]
  table[row sep=crcr]{%
16	33.51\\
32	190.17\\
64	237.09\\
128	265.28\\
256	335.2\\
};
\addlegendentry{$L = 2$}

\addplot [color=mycolor2, mark=x, mark options={solid, mycolor2}]
  table[row sep=crcr]{%
16	33.51\\
32	184.53\\
64	225.57\\
128	238.4\\
256	246.02\\
512	252.7975\\
};
\addlegendentry{$L = \log_2(n_x)-2$, retain 1}

\addplot [color=mycolor3, mark=o, mark options={solid, mycolor3}]
  table[row sep=crcr]{%
16	33.51\\
32	177.55\\
64	233.25\\
128	248.64\\
256	258.22\\
512	266.1275\\
};
\addlegendentry{$L = \log_2(n_x)-2$, retain 4}

\addplot [color=mycolor4, mark=diamond, mark options={solid, mycolor4}]
  table[row sep=crcr]{%
16	33.51\\
32	188.58\\
64	251.17\\
128	308.8\\
256	405.16\\
};
\addlegendentry{$L = \log_2(n_x)-2$, retain all}

\end{axis}
\end{tikzpicture}%
  \end{center}
  \caption{Memory usage of the preconditioner per core for the Stokes problem on a grid of size $n_x \times n_x \times n_x$, while keeping the number of levels at $L=2$, and when increasing the number of levels by 1 when $n_x$ is doubled.} \label{fig:smilu-memory}
\end{figure}
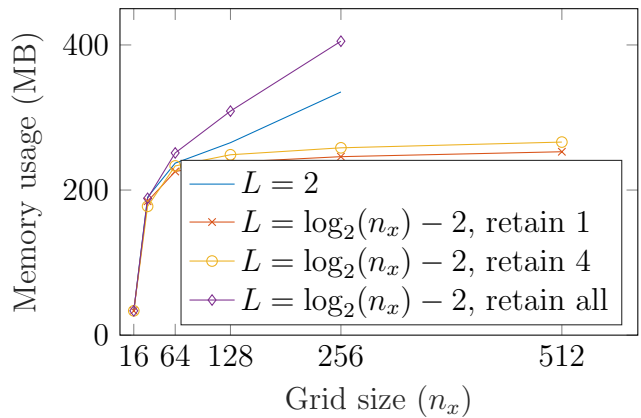

In \figref{smilu-memory}, we see the average memory usage of the preconditioner per core.
Here we again observe that the 2 level case, and the case where we retain all separator nodes after level 3, perform poorly.
The cases where we retain 1 or 4 separator nodes per separator group after level 2 show similar behavior in terms of memory usage, and the memory usage becomes constant as expected.

\subsection{Strong scalability}
In this section, we look at a problem of size $n_x=128$, with 6 levels and retaining only one node per separator group.
We use 1 to 128 cores with steps of a factor of 2 for all cases except the application time, where we use 2 to 128 cores.
The reason we do not look at 1 core for this timing is that this configuration caused a memory allocation error in the iterative solver (Belos).

\begin{figure}[!ht]
  \begin{center}
%
%
\definecolor{mycolor1}{rgb}{0.00000,0.44700,0.74100}%
\begin{tikzpicture}[%
scale=\figureScale
]

\begin{axis}[%
width=0.951\figureWidth,
height=\figureHeight,
at={(0\figureWidth,0\figureHeight)},
scale only axis,
xmin=0,
xmax=140,
xtick={1,8,16,64,128,256,512},
xlabel style={font=\color{white!15!black}},
xlabel={Number of cores ($n_p$)},
ymin=11.5,
ymax=15.5,
ylabel style={font=\color{white!15!black}},
ylabel={Memory usage (GB)},
axis background/.style={fill=white}
]
\addplot [color=mycolor1, forget plot]
  table[row sep=crcr]{%
1	11.9\\
2	14.08\\
4	14.16\\
8	14.22\\
16	14.41\\
32	14.69\\
64	14.9\\
128	15.1\\
};
\end{axis}
\end{tikzpicture}%
  \end{center}
  \caption{Total memory usage of the preconditioner for the Stokes problem on a grid of size $128^3$, with 1 to 128 cores.} \label{fig:smilu-memory-strong}
\end{figure}
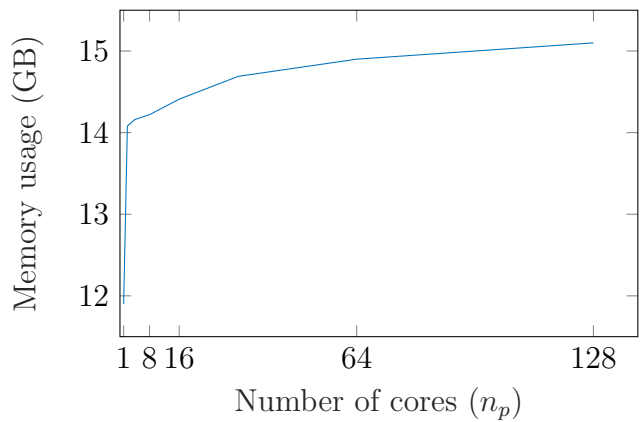

We first look at the total memory usage in \figref{smilu-memory-strong}.
We observe that there is a large difference between the memory usage on one core, and the memory usage on two cores.
This is because especially Epetra uses different implementations of its data structures for serial and parallel computations.

We also observe that the total memory usage increases slightly when using more cores.
This can be explained by the overlap that exists between different processes.
Since the difference in communication between a cubical domain and a parallelepipedal domain is a constant factor of 
$\frac{3}{4}$, we may instead look at a cubical domain to explain this behavior.
If we have a subdomain of size $s_x\times s_y \times s_z$, then communication is required for $2s_xs_y + 2s_xs_z + 2s_ys_z$ grid cells.
If the subdomain is split in the $x$-direction, communication for $2s_xs_y + 2s_xs_z + 4s_ys_z$ grid cells is required.
For the case  $s_x=s_y=s_z$, communication increases with a factor $\frac{4}{3}$ when doubling the number of cores.
This explains the slight increase in memory usage we measure.

\begin{figure}[!ht]
  \begin{center}
%
%
\definecolor{mycolor1}{rgb}{0.00000,0.44700,0.74100}%
\definecolor{mycolor2}{rgb}{0.85000,0.32500,0.09800}%
\definecolor{mycolor3}{rgb}{0.92900,0.69400,0.12500}%
\begin{tikzpicture}[%
scale=\figureScale
]

\begin{axis}[%
width=0.951\figureWidth,
height=\figureHeight,
at={(0\figureWidth,0\figureHeight)},
scale only axis,
xmin=0,
xmax=7,
xtick={0,1,2,3,4,5,6,7,8},
xticklabels={{1},{2},{4},{8},{16},{32},{64},{128},{256}},
xlabel style={font=\color{white!15!black}},
xlabel={Number of cores ($n_p$)},
ymin=0,
ymax=7,
ytick={0,1,2,3,4,5,6,7,8},
yticklabels={{1},{2},{4},{8},{16},{32},{64},{128},{256}},
ylabel style={font=\color{white!15!black}},
ylabel={Speedup},
axis background/.style={fill=white},
legend style={at={(0.03,0.97)}, anchor=north west, legend cell align=left, align=left, draw=white!15!black}
]
\addplot [color=mycolor1]
  table[row sep=crcr]{%
0	0\\
1	1\\
2	2\\
3	3\\
4	4\\
5	5\\
6	6\\
7	7\\
};
\addlegendentry{Ideal}

\addplot [color=mycolor2, mark=x, mark options={solid, mycolor2}]
  table[row sep=crcr]{%
0	0\\
1	1.14579440859251\\
2	1.88421347786198\\
3	2.74218184971564\\
4	3.68723249683006\\
5	4.61104873436546\\
6	5.56758681112997\\
7	6.51660288241116\\
};
\addlegendentry{Compute}

\addplot [color=mycolor3, mark=o, mark options={solid, mycolor3}]
  table[row sep=crcr]{%
1	1\\
2	1.97034919346174\\
3	2.85713129698795\\
4	3.83903750342608\\
5	4.78914105917907\\
6	5.62939357483162\\
7	6.61121574038262\\
};
\addlegendentry{Apply}

\end{axis}
\end{tikzpicture}%
  \end{center}
  \caption{Speedup of computation and application of the preconditioner for the Stokes problem on a grid of size $128^3$, with 2 to 128 cores.} \label{fig:smilu-timing-strong}
\end{figure}
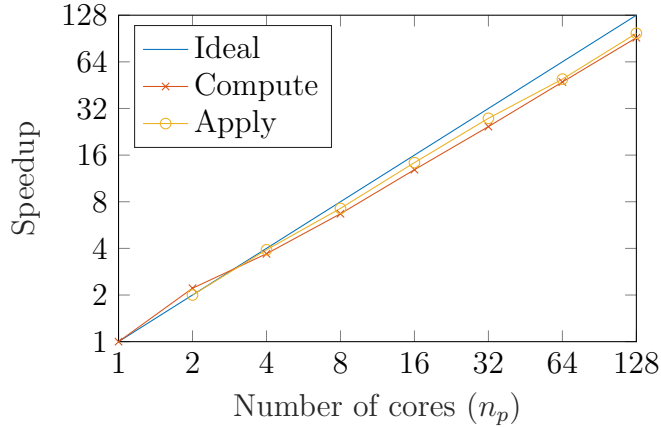

In \figref{smilu-timing-strong} we look at the speedup when using more cores.
Ideally, we would see that using twice the number of cores would mean half of the computational time.
We plotted this ideal line for reference.
We see that the speedup of both the computation and application of the preconditioner is very close to this ideal line.
There is one outlier at 2 cores which is above the ideal line, which may be explained by system variability.
Other than that the behavior is as expected.
\subsection{Lid-driven cavity}
In the previous sections we determined the weak and strong scalability properties of the preconditioner, which means that we can now continue with the robustness of the solver on the lid-driven cavity problem with increasing Reynolds numbers.
We perform a continuation with steps of $\Rey=100$ starting from the solution of the Stokes problem.
We show the results of the first iteration of Newton at $\Rey=500$ and $\Rey=2000$.
The reason we go up to $\Rey=2000$ is that a Hopf bifurcation is located between $\Rey=1900$ and 
$\Rey=2000$\citep{feldman:10,liberzon:11}, which is of interest because it changes the qualitative behavior of the solution and makes the Jacobian indefinite.

\begin{figure}[!ht]
  \begin{center}
%
%
\definecolor{mycolor1}{rgb}{0.00000,0.44700,0.74100}%
\definecolor{mycolor2}{rgb}{0.85000,0.32500,0.09800}%
\definecolor{mycolor3}{rgb}{0.92900,0.69400,0.12500}%
\begin{tikzpicture}[%
scale=\figureScale
]

\begin{axis}[%
width=0.951\figureWidth,
height=\figureHeight,
at={(0\figureWidth,0\figureHeight)},
scale only axis,
xmin=0,
xmax=600,
xtick={ 16,  64, 128, 256, 512},
xlabel style={font=\color{white!15!black}},
xlabel={Grid size ($n_x$)},
ymin=0,
ymax=2500,
ylabel style={font=\color{white!15!black}},
ylabel={Number of iterations},
axis background/.style={fill=white},
legend style={at={(0.03,0.97)}, anchor=north west, legend cell align=left, align=left, draw=white!15!black}
]
\addplot [color=mycolor1]
  table[row sep=crcr]{%
16	171\\
32	340\\
64	508\\
128	603\\
256	657\\
};
\addlegendentry{$L = 2$}

\addplot [color=mycolor2, mark=x, mark options={solid, mycolor2}]
  table[row sep=crcr]{%
16	171\\
32	364\\
64	619\\
128	1031\\
256	1546\\
512	2249\\
};
\addlegendentry{$L = \log_2(n_x)-2$, retain 1}

\addplot [color=mycolor3, mark=o, mark options={solid, mycolor3}]
  table[row sep=crcr]{%
16	171\\
32	364\\
64	619\\
128	1017\\
256	1464\\
512	1942\\
};
\addlegendentry{$L = \log_2(n_x)-2$, retain 4}

\end{axis}
\end{tikzpicture}%
  \end{center}
  \caption{Number of iterations of GMRES for the lid-driven cavity problem at $\Rey=500$ on a grid of size $n_x \times n_x \times n_x$, while keeping the number of levels at $L=2$, and when increasing the number of levels by 1 when $n_x$ is doubled.
    A relative residual of $10^{-8}$ was used as convergence tolerance.
    As initial guess we used the solution at $\Rey=400$.} \label{fig:ldc-its-500}
\end{figure}
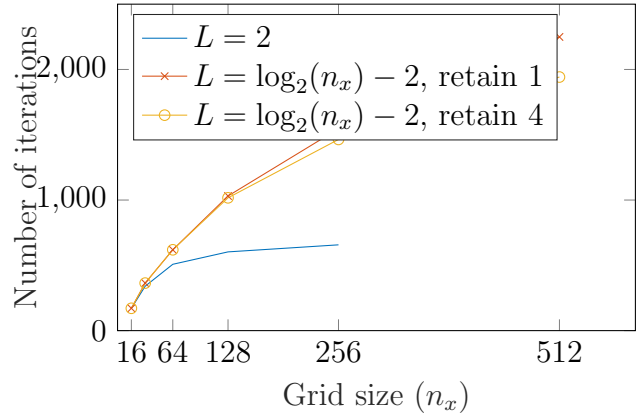

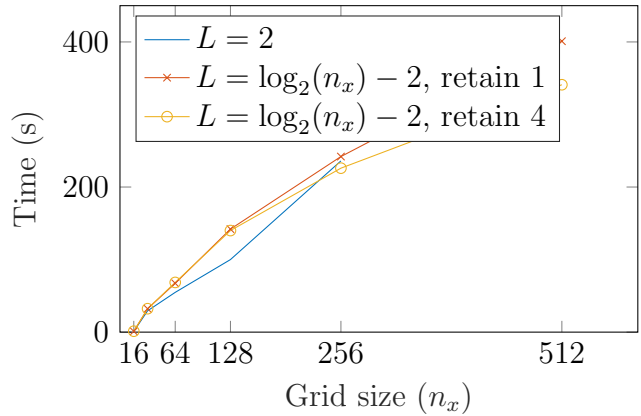
\begin{figure}[!ht]
  \begin{center}
%
%
\definecolor{mycolor1}{rgb}{0.00000,0.44700,0.74100}%
\definecolor{mycolor2}{rgb}{0.85000,0.32500,0.09800}%
\definecolor{mycolor3}{rgb}{0.92900,0.69400,0.12500}%
\begin{tikzpicture}[%
scale=\figureScale
]

\begin{axis}[%
width=0.951\figureWidth,
height=\figureHeight,
at={(0\figureWidth,0\figureHeight)},
scale only axis,
xmin=0,
xmax=600,
xtick={ 16,  64, 128, 256, 512},
xlabel style={font=\color{white!15!black}},
xlabel={Grid size ($n_x$)},
ymin=0,
ymax=450,
ylabel style={font=\color{white!15!black}},
ylabel={Time (s)},
axis background/.style={fill=white},
legend style={at={(0.03,0.97)}, anchor=north west, legend cell align=left, align=left, draw=white!15!black}
]
\addplot [color=mycolor1]
  table[row sep=crcr]{%
16	1.07\\
32	29.4\\
64	54.8\\
128	100\\
256	236\\
};
\addlegendentry{$L = 2$}

\addplot [color=mycolor2, mark=x, mark options={solid, mycolor2}]
  table[row sep=crcr]{%
16	1.07\\
32	32.2\\
64	67.3\\
128	142\\
256	242\\
512	401\\
};
\addlegendentry{$L = \log_2(n_x)-2$, retain 1}

\addplot [color=mycolor3, mark=o, mark options={solid, mycolor3}]
  table[row sep=crcr]{%
16	1.27\\
32	32.3\\
64	68.5\\
128	140\\
256	226\\
512	341\\
};
\addlegendentry{$L = \log_2(n_x)-2$, retain 4}

\end{axis}
\end{tikzpicture}%
  \end{center}
  \caption{Time for GMRES to converge for the lid-driven cavity problem at $\Rey=500$ on a grid of size $n_x \times n_x \times n_x$, while keeping the number of levels at $L=2$, and when increasing the number of levels by 1 when $n_x$ is doubled.
    A relative residual of $10^{-8}$ was used as convergence tolerance.
    As initial guess we used the solution at $\Rey=400$.} \label{fig:ldc-time-500}
\end{figure}

In \figref{ldc-its-500} and \figref{ldc-time-500}, we show the results at Reynolds number 500, which show good correspondence with the results on the Stokes problem.
The main difference is that more iterations are needed, since higher Reynolds numbers constitute harder problems.
We see that the number of iterations that is required when keeping the number of levels constant actually decreases, which is due to the grid refinement having a positive effect on the mesh P\'eclet number, i.e.\ the coefficients of the diffusion increase with respect to those of the convection.

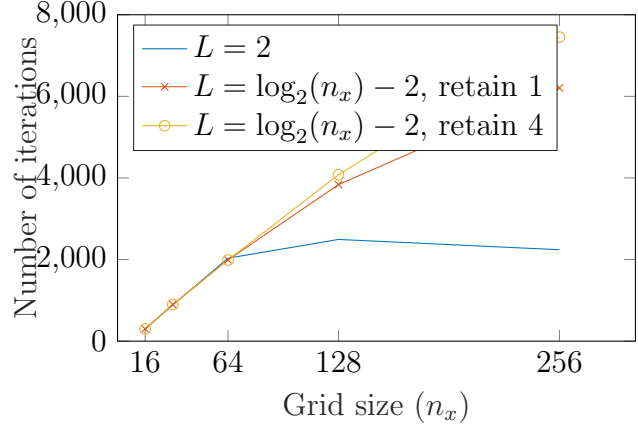
\begin{figure}[!ht]
  \begin{center}
%
%
\definecolor{mycolor1}{rgb}{0.00000,0.44700,0.74100}%
\definecolor{mycolor2}{rgb}{0.85000,0.32500,0.09800}%
\definecolor{mycolor3}{rgb}{0.92900,0.69400,0.12500}%
\begin{tikzpicture}[%
scale=\figureScale
]

\begin{axis}[%
width=0.951\figureWidth,
height=\figureHeight,
at={(0\figureWidth,0\figureHeight)},
scale only axis,
xmin=0,
xmax=300,
xtick={ 16,  64, 128, 256, 512},
xlabel style={font=\color{white!15!black}},
xlabel={Grid size ($n_x$)},
ymin=0,
ymax=8000,
ylabel style={font=\color{white!15!black}},
ylabel={Number of iterations},
axis background/.style={fill=white},
legend style={at={(0.03,0.97)}, anchor=north west, legend cell align=left, align=left, draw=white!15!black}
]
\addplot [color=mycolor1]
  table[row sep=crcr]{%
16	300\\
32	886\\
64	2033\\
128	2493\\
256	2241\\
};
\addlegendentry{$L = 2$}

\addplot [color=mycolor2, mark=x, mark options={solid, mycolor2}]
  table[row sep=crcr]{%
16	300\\
32	898\\
64	1990\\
128	3837\\
256	6206\\
};
\addlegendentry{$L = \log_2(n_x)-2$, retain 1}

\addplot [color=mycolor3, mark=o, mark options={solid, mycolor3}]
  table[row sep=crcr]{%
16	300\\
32	898\\
64	1990\\
128	4079\\
256	7451\\
};
\addlegendentry{$L = \log_2(n_x)-2$, retain 4}

\end{axis}
\end{tikzpicture}%
  \end{center}
  \caption{Number of iterations of GMRES for the lid-driven cavity problem at $\Rey=2000$ on a grid of size $n_x \times n_x \times n_x$, while keeping the number of levels at $L=2$, and when increasing the number of levels by 1 when $n_x$ is doubled.
    A relative residual of $10^{-8}$ was used as convergence tolerance.
    As initial guess we used the solution at $\Rey=1900$.} \label{fig:ldc-its-2000}
\end{figure}

\begin{figure}[!ht]
  \begin{center}
%
%
\definecolor{mycolor1}{rgb}{0.00000,0.44700,0.74100}%
\definecolor{mycolor2}{rgb}{0.85000,0.32500,0.09800}%
\definecolor{mycolor3}{rgb}{0.92900,0.69400,0.12500}%
\begin{tikzpicture}[%
scale=\figureScale
]

\begin{axis}[%
width=0.951\figureWidth,
height=\figureHeight,
at={(0\figureWidth,0\figureHeight)},
scale only axis,
xmin=0,
xmax=300,
xtick={ 16,  64, 128, 256, 512},
xlabel style={font=\color{white!15!black}},
xlabel={Grid size ($n_x$)},
ymin=0,
ymax=1500,
ylabel style={font=\color{white!15!black}},
ylabel={Time (s)},
axis background/.style={fill=white},
legend style={at={(0.03,0.97)}, anchor=north west, legend cell align=left, align=left, draw=white!15!black}
]
\addplot [color=mycolor1]
  table[row sep=crcr]{%
16	2.5\\
32	81\\
64	228\\
128	409\\
256	860\\
};
\addlegendentry{$L = 2$}

\addplot [color=mycolor2, mark=x, mark options={solid, mycolor2}]
  table[row sep=crcr]{%
16	2.5\\
32	80\\
64	212\\
128	550\\
256	980\\
};
\addlegendentry{$L = \log_2(n_x)-2$, retain 1}

\addplot [color=mycolor3, mark=o, mark options={solid, mycolor3}]
  table[row sep=crcr]{%
16	3.1\\
32	83\\
64	220\\
128	580\\
256	1230\\
};
\addlegendentry{$L = \log_2(n_x)-2$, retain 4}

\end{axis}
\end{tikzpicture}%
  \end{center}
  \caption{Time for GMRES to converge for the lid-driven cavity problem at $\Rey=2000$ on a grid of size $n_x \times n_x \times n_x$, while keeping the number of levels at $L=2$, and when increasing the number of levels by 1 when $n_x$ is doubled.
    A relative residual of $10^{-8}$ was used as convergence tolerance.
    As initial guess we used the solution at $\Rey=1900$.} \label{fig:ldc-time-2000}
\end{figure}
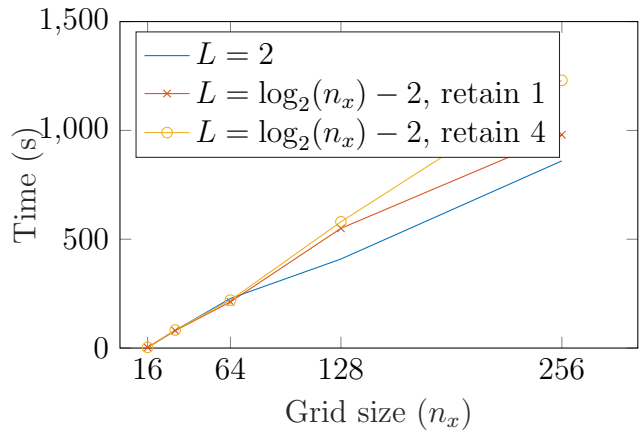

The results at Reynolds number 2000 are shown in \figref{ldc-its-2000} and \figref{ldc-time-2000}.
We again observe that the number of iterations is much larger.
What is odd, however, is that retaining more nodes now actually gives worse convergence.
This may be because we use GMRES(250) instead of GMRES due to memory limitations, and therefore do not preserve the convergence properties of GMRES.
This effect is more prevalent for this problem because of the large number of iterations that is required.
We did confirm that for the first 250 iterations, retaining 4 nodes after level 2 gives rise to better convergence, as we expected.
For Reynolds number 2000, we do not show results with a grid of size $512^3$, because the method does not converge within $10\,000$ GMRES iterations, which we set as the maximum.
It does converge, however, and extrapolating the results, we expect that around 15000 iterations are required to meet the tolerance.
Alternatively, instead of increasing the number of levels with the grid size, we could also have used 6 levels as we did for a grid of size $256^3$, in which case we would have expected it to actually use fewer iterations due to the effect of grid refinement on the mesh P\'eclet number as explained earlier.

In \tabref{teko} we show results for Reynolds number 500 using the LSC-SIMPLEC block preconditioner implemented in Teko \citep{cyr:12}.
We also tried preconditioners from other packages, as well as other preconditioning strategies implemented in the Teko package, but those unfortunately converged even slower or did not yield convergence at all.
The stopping criterion of the linear solver is $10^{-8}$, as before.

\begin{table}[ht]
\begin{center}
\begin{tabular}{rrrrrrr}
          $n_p$    &    $n_x$   &       $L$  &     its  &      $t_c$ (s)    &       $t_s$ (s)\\\hline
        1 &        16 &         2 &       142 &      0.12 &      0.77 \\
        1 &        32 &         2 &       187 &      9.88 &     18.80 \\
        8 &        64 &         3 &       245 &   1511.12 &    313.00
\end{tabular}
\end{center}
\caption{Performance of Teko with LSC-SIMPLEC preconditioner for the lid-driven cavity problem at $\Rey = 500$ on a grid of size $n_x \times n_x \times n_x$.
  Here $n_p$ is the number of cores, $L$ is the number of levels, its is the number of iterations, $t_c$ is the time to compute the preconditioner and $t_s$ is the time for solving the linear system.
  \label{tab:teko}}
\end{table}

Compared to our method, we see that Teko has much more difficulty with the grid refinement, leading to a huge computational cost at a grid size of only $64^3$.
A crude computation shows that per grid point the method becomes about 65 times more expensive per iteration.
This must be attributed to slow convergence of algebraic multigrid on the subblocks.
One of these blocks is the $L+N$ block from \eqref{ns_spp}, with $N$ indefinite, and this is something that is difficult for a standard AMG method.
We choose the number of levels to be 2 for grid sizes $16^3$ and $32^3$, and 3 levels for $64^3$ since these seem to give the optimal results.
For Reynolds number 2000, we do not observe convergence past the $16^3$ grid.
\section{Summary and Discussion}\label{sec:smilu-conc}
We presented a robust method for solving the steady, incompressible Navier--Stokes equations, which makes use of parallelepiped shaped overlapping subdomains.
The interiors of these overlapping subdomains can be eliminated in parallel.
On the interfaces of the subdomains, Householder transformations are applied to decouple all but one velocity node from the pressure nodes, after which all decoupled nodes can also be eliminated in parallel.
The key to the multilevel approach is the resulting reduced Schur complement matrix, which has the same structure as the original matrix.
We can therefore recursively apply our method to this matrix.
The shape of the subdomains makes sure that pressure nodes are not isolated in the factorization process, which would result in a singular Schur complement matrix.

Our weak scalability experiments show that if the number of levels of the multilevel approach is kept constant while increasing the problem size, grid independent convergence is obtained.
We also show that by increasing the number of levels and processors as the problem size increases, we only require a small amount of additional time and memory for both the factorization and solution process.
Moreover, by increasing the number of nodes that is retained per separator after level 2, we can further decrease the time that is required to solve the system.

Our strong scalability experiments show that the time that is required to both compute and apply the preconditioner scales very well.
The same holds for the memory usage, which behaves as expected.

We also showed that the preconditioner leads to convergence of GMRES for the lid-driven cavity problem at high Reynolds numbers, where other methods, such as the LSC-SIMPLEC preconditioner that is implemented in Teko, fail to do so.

This leads us to conclude that we implemented a robust solution method for the Navier--Stokes equations on staggered grids which shows good parallel scalability.
In this paper, we only showed results for Arakawa C-grids.
We have already succeeded to apply the idea to discretization on other staggered grids, e.g.\ Arakawa B-grids and mixtures of B- and C-grids as used in oceanography.
This will allow us to use it as preconditioner in ocean-climate models, which is our application goal.
The approach is, however, not limited to structured grids.
The essence of the method is that we can iterate in the divergence-free space in a cheap way.
As far as we can see this is the case when we have discrete conservation of mass.
The C-grid discretizaton leads to one mass conservation law per subdomain (the B-grid has two).
By uniting two neighbouring subdomains, we can find, by a simple combination of the two laws, a single new one for the united subdomains.
This is also possible for finite volume discretizations on unstructured grids.
For this reason, we think that the method can be applied to such discretizations in general and to some (discontinuous) Galerkin methods that also have the discrete conservation property.
In the future we may generalize the approach to unstructured grid discretizations for which we will need graph based partitioning methods.

\subsection*{Acknowledgments}
This work is part of the Mathematics of Planet Earth research program with project number 657.014.007, which is financed by the Netherlands Organization for Scientific Research (NWO) and the SMCM project of the Netherlands eScience Center (NLeSC) with project number 027.017.G02 (SB and FW).
Access to the SuperMUC-NG system during the `friendly user phase' was kindly granted by LRZ Munich via a test account.
We would like to thank Eric Cyr for his help during our attempts of getting the solvers in Teko to work as expected.

\bibliographystyle{abbrvnatdoi}
\bibliography{smilu-paper}

\end{document}